\documentclass[reqno]{amsart}
\usepackage{amssymb,latexsym}
\usepackage{amsfonts,mathrsfs}
\usepackage{amsmath}
\usepackage{txfonts}
\usepackage[margin=1.4in]{geometry}

\newtheorem{theorem}{Theorem}[section]
\newtheorem{lemma}[theorem]{Lemma}
\newtheorem{proposition}[theorem]{Proposition}
\newtheorem{corollary}[theorem]{Corollary}
\newtheorem{definition}[theorem]{Definition}

\theoremstyle{remark}
\newtheorem{remark}{Remark}[section]

\newcommand{\Vol}{\mathrm{Vol}}

\begin{document}
\title{ The  Prescribed $Q$-Curvature Flow for Arbitrary Even Dimension in a Critical Case}
 \author{Yuchen Bi} \address{ Institute of Mathematics, Academy of Mathematics and Systems
  Science, University of Chinese Academy of Sciences, Beijing, 100190, P. R.
China} \email{biyuchen15@mails.ucas.ac.cn}
\author{Jiayu Li}
\address{School of Mathematical Sciences, University of Science and Technology of China,
Hefei, 230026, P. R. China} \email{jiayuli@ustc.edu.cn}
\thanks {The research was supported by the National Key Research and Development Project SQ2020YFA070080 and NSFC No 11721101 }
\thanks{ }

\begin{abstract}
In this paper, we study the prescribed $Q$-curvature flow equation on a arbitrary even dimensional closed Riemannian manifold $(M,g)$, which was introduced by S. Brendle in \cite{B2003}, where he proved the flow  exists for long time and converges at infinity if  the GJMS operator is weakly positive with trivial kernel and $\int_M Qd\mu < (n-1)!\Vol\left( S^n \right) $. In this paper we study the critical case that $\int_M Qd\mu = (n-1)!\Vol\left( S^n \right)$,  we will prove the convergence of the flow under some geometric hypothesis.  In particular, this gives a new proof of Li-Li-Liu's existence result in \cite{LLL2012} in dimensiona 4 and extend the work of Li-Zhu \cite{LZ2019} in dimension 2 to general even dimensions. In the proof, we give a explicit  expression of the limit of the corresponding  energy functional  when the blow up occurs.
\end{abstract}

\maketitle

\numberwithin{equation}{section}
\section{Introduction}

Consider a compact four dimensional manifold $(M,g)$. The $Q$-curvature is defined by
$$Q=-\frac{1}{6}(\Delta R +R^2-3|Ric|^2) ,$$
where $R$ denotes the scalar curvature and $Ric$ denotes the Ricci tensor. Under a conformal change of the metric $\tilde{g}= e^{2u}g$, the following equation holds:
\begin{equation}\label{trans}
  \tilde{Q}=e^{-4u}(Pu+Q) ,
\end{equation}
where $P$ denotes the Paneitz operator defined by
$$Pu=\Delta^2u+div\left[\left(\frac{2}{3}Rg-2Ric\right)\nabla u\right].$$
The Paneitz operator is conformally covariant in the sense that if $\tilde{g}= e^{2u}g$, then the Paneitz operator with respect to $\tilde{g}$ is given by $\tilde{P}=e^{-4u}P$. The Gauss-Bonnet-Chern theorem asserts that
$$\int_M Qd\mu+\frac{1}{4}\int_M|W|^2d\mu=8\pi^2\chi(M),$$
where $\chi(M)$ denotes the Euler characteristic of $M$, the letter $W$ denotes the Weyl curvature tensor of $g$ and $d\mu$ denotes the volume element  of $g$.  Since the Weyl curvature tensor is conformally invariant, it follows that
\begin{equation}\label{intinv}
  \int_M Qd\mu=\int_M \tilde{Q}d\tilde{\mu}.
\end{equation}
We note here that the invariance (\ref{intinv}) is also a direct consequence of equation (\ref{trans}).

So $Q$ and $P$ are the four dimensional analogues  to the Gaussian curvature $K$ and the Laplace–Beltrami operator $-\Delta$ in the two dimensional case. We denote $Q=K$ and $P= -\Delta$ in two dimensional for simplicity.

In \cite{GJMS1992}, it  was  introduced  a family  of  conformally  covariant  differential  operators on general even dimensional manifolds,  whose leading term is $\left( -\Delta \right) ^{\frac{n}{2}}$, were $n$ is the dimension of the manifold. These operators are usually referred to as the GJMS operators. We denote them also by $P$. In \cite{B1993}, the corresponding $Q$-curvature was defined.

Furthermore as for the Laplace–Beltrami operator on compact surfaces and the Paneitz operator on compact four-dimensional manifold, for every compact $n$-dimensional manifolds with $n$ even, we have that after a conformal change of metric $\tilde{g} = e^{2u}g$.
\[
	\tilde{P} = e^{-nu}P\qquad\qquad \tilde{Q}=e^{-nu}(Pu+Q)
.\]

An important problem in conformal geometry is the construction of a conformal metric $\tilde{g}= e^{2u}g$ for which the  $Q$-curvature $\tilde{Q}$ equals to a constant multiple of a prescribed function $f$. This is related to solve the equation
\begin{equation}\label{elliptic equation}
  Pu+Q=lfe^{nu}.
\end{equation}
From the equation (\ref{intinv}), we know that $l=\frac{\bar{Q}}{\bar{f}}$ , where
$$\bar{Q}=\int_MQd\mu;\qquad\bar{f}=\int_M fd\tilde{\mu}=\int_Mfe^{nu}d\mu.$$
The equation (\ref{elliptic equation}) is the Euler-Lagrange equation of the functional
\begin{equation}\label{functional}
	E[u]=\frac{n}{2}\int_MuPud\mu+n\int_M Q ud\mu-\bar{Q}\log\left(\int_Mfe^{nu}d\mu\right).
\end{equation}

S.-Y. Chang and P. C. Yang \cite{CY1995} first studied the equation (\ref{elliptic equation}) in the case $n=4$ by minimizing the functional $E$. They constructed conformal metrics of constant $Q$-curvature when the Paneitz operator $P$ is weakly positive with trivial kernel and the total $Q$-curvature satisfies $\bar{Q}< 16\pi^2$. The Paneitz operator  $P$ is weakly positive means that for all $u\in C^{\infty}(M)$, $\int_MuPud\mu\geq0$, and $P$ has trivial kernel means that its kernel consists only of constant functions. In view of the result of M. Gursky \cite{G1999}, $P$ is positive with trivial kernel whenever $\bar{Q}>0$ and the Yamabe constant of $(M,g)$ is positive. The key point in S.-Y. Chang and P. C. Yang's proof is the following: Using Adams-Fontana inequality (see Proposition \ref{AF Proposition} below), supposing $P$ is weakly positive with trivial kernel and $\bar{Q}< 16\pi^2$, they deduced that the functional $E$ is bounded from below and coercive. Then the solutions of equation (\ref{elliptic equation}) can be  found as global minima of $E$.

Later, using the heat flow method, S. Brendle \cite{B2003} extended S.-Y. Chang and P. C. Yang' result to higher dimensions and to $f>0$ which need not  be constants. We now sketch S. Brendle's proof in the four dimension case. First he introduced the flow equation:
\begin{equation}\label{parabolic equation}
  \begin{cases}
    \frac{\partial u}{\partial t}=-e^{-nu}(Pu+Q)+\frac{\bar{Q}}{\bar{f}}f\\
    u(x,0)=u_0(x)\in C^{\infty}(M),
  \end{cases}
\end{equation}
where $\bar{Q}=\int_M Qd\mu$ and $\bar{f}(t)=\int_Mfe^{nu(t)}d\mu$. It is the negative gradient flow of the functional $E$ with respect to the inner product $\int_M(-,-)e^{nu(t)}d\mu$. One advantage of this flow is that it keeps the total volume constant, that is, $\int_M e^{nu(t)}d\mu\equiv constant$ along this flow. Then he derived various differential inequalities from the Gagliardo-Nirenberg inequality. Combining these inequalities and a stability argument, Brendle obtained that (notice that we state  the result of  S.Brendle  in a form more convenient for us):
\begin{theorem}\label{brendle main}
	Assume that the GJMS operator $P$ is weakly positive with trivial kernel, and $u(t)$ is the solution of the equation (\ref{parabolic equation}). Then:
	\\
(a) If we also assume that $\sup_{t\in[0,T_0]}||u(t)||_{W^{\frac{n}{2},2}(M)}\leq C(g,f,T_0)< +\infty$ for any  $T_0\in[0,T)$, where $T$ is the maximal existence time.  Then $T=+\infty$, that is, $u(t)$ exists long time.\\
(b) If we also assume that $\sup_{t\in[0,T_0]}||u(t)||_{W^{\frac{n}{2},2}(M)}\leq C(g,f)< +\infty$ for any  $T_0\in[0,T)$, where $T$ is the maximal existence time and $C(g,f)$ is a constant which is independent of $T_0$. Then $u(t)$ exists long time and  $u(t)$ converges smoothly to a smooth function $u_{\infty}$ satisfying the equation (\ref{elliptic equation}).  \end{theorem}

Then by using the coerciveness,  he concluded that the condition of part $(b)$ in above proposition is satisfied when $\int_M Qd\mu< (n-1)!\gamma_n$. Here we use $\gamma_n$ to denote the volume of the standard round $n$-sphere. As a consequence, in the case that $\int_M Qd\mu<(n-1)!\gamma_n$, he proved long time existence and convergence of the flow equation (\ref{parabolic equation}) and the limit at infinity is a solution of the equation (\ref{elliptic equation}). In this way, S.Brendle generalized S.-Y. Chang and P. C. Yang' existence result.  For more information about the flow equation (\ref{parabolic equation}), one can refer to \cite{B2003}, \cite{B2006}, \cite{MS2006} and \cite{FR2018}.
In the critical case $\bar{Q}=(n-1)!\gamma_n$, the study of the equation (\ref{elliptic equation}) becomes more subtle. The main difficulty is the appearance of the so-called bubbling phenomena due to the concentration of the volume of the conformal metric. Note that when $(M, g) = (S^n, g_{round})$, the existence of (\ref{elliptic equation}) is an analogue of the well-known Nirenberg's problem. It has been studied by S. Brendle \cite{B2006}  and A. Malchiodi \& M. Struwe \cite{MS2006} by using the flow methods.

One can also consider the equation (\ref{elliptic equation}) replacing the $Q$ term by an arbitrary smooth function, which is called the mean field type equation. See \cite{CL2002}, \cite{CL2003}, \cite{M2009},\cite{N2017} for more information.

In this paper we study the flow equation (\ref{parabolic equation})
in the critical case for general Riemannian manifolds $M$. Our results can apply to the mean field type equations. Notice that the $n=2$ case have been studied by the second author and C. Zhu in \cite{LZ2019}.

Let  $G_p$ be the unique solution of the following equation
\begin{equation}\label{Green}
    \begin{cases}
	    PG_p+ Q= (n-1)!\gamma_n\delta_p\\
      \int_M G_pd\mu=0,
    \end{cases}
  \end{equation}
  that is, $G_p$ is the Green's function with a pole at $p$.

We have
\begin{theorem}\label{two alternatives}
	Let $(M,g)$ be a closed Riemannian manifold of even dimension $n\geq 2$, with $\bar{Q}=(n-1)!\gamma_n$. Suppose the GJMS operator $P$ is weakly  positive with trivial kernel and $f>0$.  Assume $u(t)$ is the solution of the parabolic equation:
  \begin{equation}
    \begin{cases}
      \frac{\partial u}{\partial t}=-e^{-nu}(Pu+Q)+\frac{\bar{Q}}{\bar{f}}f\\
      u(x,0)=u_0(x)\in C^{\infty}(M),
    \end{cases}
  \end{equation}
  where $\bar{Q}=\int_M Qd\mu$ and $\bar{f}(t)=\int_Mfe^{nu(t)}d\mu$. Then $u(t)$ exists for long time and one of the following alternatives holds:\\

  (a) either $u(t)$ converges smoothly to a smooth function $u_{\infty}$ satisfying the equation (\ref{elliptic equation}),

  (b) or
  \begin{equation}\label{lower energy bound}
    \begin{split}
      \lim_{t\to+\infty}E[u(t)]& =	-(n-1)!\gamma_n\log\frac{\gamma_n f(p)}{2^n}-\frac{n!\gamma_n}{2}\left( H_p(p)+\frac{C_{\log}}{C_0} \right)+\frac{n}{2}\int_M QG_pd\mu\\
			       &=: \Lambda(g,f,p)
    \end{split}
  \end{equation}
  for some $p\in M$, where  $\frac{C_{\log}}{C_0}$ is a uniform constant, $G_p$ is the Green's function with a pole at $p$, and the $H_p$  is the regular part of the Green's function as in Proposition \ref{Green expan}.
\end{theorem}

The new feature of our results is that we can show $\lim_{t\to\infty}E[u(t)]=\Lambda\left( g,f,p \right)$ for the alternative (b). Even in the two dimension, only  $\lim_{t\to\infty}E[u(t)]\geq \Lambda\left( g,f,p \right)$ has been known, see J. Li and C. Zhu in \cite{LZ2019}.

We can show the convergence of this flow with arbitrary initial data under some geometric conditions in the critical case. To the authors knowledge, similar results has only obtained when $M$ is the standard $S^n$ ( c.f. \cite{B2006}, \cite{H2010} ) .

\begin{theorem}\label{large convergence theorem}
  Under the assumption of Theorem \ref{two alternatives}, if
  \begin{equation}
	  \left( \Delta\left( H_p+\frac{1}{n}\log f \right)+n\left| \nabla\left( H_p+\frac{1}{n}\log f \right)  \right|^2 -\frac{R_g}{6(n-1)}  \right)(p) >0
  \end{equation}
  for all $p\in M$, and
  \begin{equation}\label{constant energy}
	  \Lambda\left( g, f, p \right) \equiv \Lambda
  \end{equation}
   for some constant $\Lambda$ as a function of $p$,   then for
  arbitrary initial data $u_0\in C^{\infty}(M)$ ,  $u(t)$ exists for long time and converges smoothly to a smooth function $u_{\infty}$ satisfying the equation (\ref{elliptic equation}).
\end{theorem}
\begin{remark}
	When $M$ is the two dimensional flat torus $T^2$, $Q$ is a constant, $\bar{Q}=4\pi$ and $f\equiv 1$, using (\ref{2dim}), it is easy to check the conditions in Theorem \ref{large convergence theorem} are satisfied. This provide an interesting analogue of the Ricci flow on $S^2$.
\end{remark}

We also have the convergence of (\ref{parabolic equation}) with small arbitrary initial data :
\begin{theorem}\label{convergence theorem}
  Under the assumption of Theorem \ref{two alternatives}, if
  \begin{equation}\label{sufficient condition}
	  \left( \Delta\left( H_{p_0}+\frac{1}{n}\log f \right)+n\left| \nabla\left( H_{p_0}+\frac{1}{n}\log f \right)  \right|^2 -\frac{R_g}{6(n-1)}  \right)(p_0) >0
  \end{equation}
  for some $p_0$ where $\Lambda(g,f,p)$ achieves its global minimum,   then there exists an initial data $u_0\in C^{\infty}(M)$ such that $u(t)$  exists  for long time and converges smoothly to a smooth function $u_{\infty}$ satisfying the equation (\ref{elliptic equation}).
  \end{theorem}

  Theorem \ref{convergence theorem} extends the main result in \cite{LZ2019} to general even dimensions and provides a new proof of the existence result of J. Li et al. (Theorem 1.2 in \cite{LLL2012}).

\setcounter{equation}{0}
\hspace{0.4cm}

\noindent {\bf Acknowledgements}.
The authors thank Dr. Chaona Zhu for her helpful discussion. The first author is also grateful to
Dr. Jie Zhou for his warm encouragement.
\section{Preliminaries }

In this section we collect some useful preliminary facts and analyze the volume concentration phenomena of the equation (\ref{parabolic equation}).

Actually we study the following perturbation equation:
\begin{equation}\label{sequence equation}
  \begin{cases}
  Pu +Q = (lf+ h)e^{nu},
\\
\int_M lfe^{nu}d\mu = \int_M Qd\mu = (n-1)!\gamma_n ,
  \end{cases}
\end{equation}
where $l>0$ is a positive real number, $f\in C^{ \infty } (M)$ is a positive real function and $h$ is a real function with
$$E_\alpha =\int_M h^\alpha e^{nu}d\mu < \infty$$
for some $\alpha\in(1,\infty)$.

Notice that the equation (\ref{parabolic equation}) is a special case of (\ref{sequence equation}) with $h=\frac{\partial u}{\partial t}$.

We care about the behavior of $u$ as $E_\alpha$ tends to zero.

To study the equation (\ref{sequence equation}), we begin with two propositions.

The following proposition concerning solutions of (\ref{sequence equation}) is proved in \cite{M2006} by using the asymptotic behavior of the Green's function near the pole, see also \cite{DR2006},  \cite{FR2018}.
\begin{proposition}\label{weak global estimate}
	Let $(u,F)\in W^{n,1} (M)\times L^1(M)$ satisfy
  $$Pu=F$$
  with $||F||_{L^1(M)}\leq K$ for some constant $K$. Then for any $x\in M$, for any $r>0$, for any $j\in [1, n-1]\cap\mathbb{N}$ and $q\in[1,n/j)$, we have
  $$\int_{B_r(x)}|\nabla^j u|^qd\mu\leq Cr^{n-jq},$$
  where C is a positive constant depending on $K$, $M$, $q$.
\end{proposition}

Furthermore, we need the following proposition proved in \cite{M2006}:

\begin{proposition}\label{weak compactness}
	Let $(u_k, F_k)\in W^{n,1} (M)\times L^1(M)$ satisfy
  $$ Pu_k=F_k$$
  with $||F_k||_{L^1(M)}\leq K$ for some constant $K$ independent of $k$. Then:

  (a) either
  $$\int_Me^{nq(u_k-\bar{u}_k)}d\mu\leq C$$
  for some $q>1$ and some $C=C(n,K,q)>0$,

  (b) or there exists a point $x\in M$ such that for any $r>0$, we have
  $$\liminf_{k\to+\infty}\int_{B_r(x)}|F_k|d\mu\geq \frac{1}{2}(n-1)!\gamma_n.$$
\end{proposition}

\begin{remark}\label{remark weak compactness}
  Proposition \ref{weak compactness} remains valid if one replaces the metric $g$ on $M$ by a family of metrics $(g_k)_k$ depending on $k$ which is uniformly bounded in $C^m(M)$ for any $m\in\mathbb{N}$. The same result also holds if we replace $M$ by any precompact open set in a noncompact manifold and assuming all the functions with compact support in this open set.
\end{remark}

Using above two propositions, we have following estimate independent of specific form of $f$ and $h$:
\begin{proposition}\label{strong compactness}
	Let $U$ be a precompact open set in M, and $u$ be the solution of the equation (\ref{sequence equation}). There is a uniform $A\in \mathbb{R}$ such that when $E_\alpha\leq A$, we have the following:
	
	If there is some $s>0$, such that $\int_{B_s(x)} lfe^{nu}d\mu\leq \frac{1}{8}(n-1)!\gamma_n $ for every $x\in U$, then for any compact subset $K$ of $U$,  we have
	$$||u-a_{U}||_{W^{ n, \alpha } (K)}\leq C,$$
where $C=C(K, U, n, \alpha, t)$ is a constant and $a_{U}=\fint_{U}ud\mu$.
\end{proposition}

\medskip

{\bf Proof. } Choose precompact open sets $V$ and $W$ satisfying
$$K\subset W, \quad \overline{W}\subset V, \quad \overline{V}\subset U.$$
and a smooth cut-off function $\eta$ satisfying
$$\begin{cases}
  \eta(x)=1\quad x\in \overline{V};\\
  \eta(x)=0\quad x\notin U.
\end{cases}$$
We also set $w= \eta(u-a_{U})$. By Proposition \ref{weak global estimate}, the precompactness of $U$ and the boundness of $l$ one finds
$$\int_U |\nabla^ju|^qd\mu\leq C_U\quad q\in\left[1,n/j\right),\quad j\in [1, n-1]\cap\mathbb{N}$$
and hence by the Poincar\'{e} inequality (recall that $w$ has a uniform compact support and $(u-a_{U})$ has mean value $0$ on $U$) it follows that
\begin{equation}\label{soboelev bound 1}
	||w||_{W^{n-1,q}(U)}\leq C_U\quad q\in\left(1,\frac{n}{n-1}\right),
\end{equation}
and
\begin{equation}\label{soboelev bound 2}
  ||u-a_{U}||_{W^{n-1,q}(U)}\leq C_U\quad q\in\left(1,\frac{n}{n-1}\right).
\end{equation}
By equation (\ref{sequence equation}) there holds
\begin{equation}\label{equation in key proposition}
  \begin{split}
	  Pw&= \left( -\Delta\right)^{ \frac{n}{2} }  w + \mathcal{O}\left( \sum_{j=0}^{n-1}|\nabla w| \right)  \\
	    &=\eta\left( -\Delta\right)^{ \frac{n}{2} }  u + \mathcal{O}\left( \sum_{j=0}^{n-1}\left(|\nabla^j w| + |\nabla^j (u-a_{U})|\right)\right) \\
	&=\eta(lf+h)e^{nu}+ \mathcal{O}\left( \sum_{j=0}^{n-1}\left(|\nabla^j w| + |\nabla^j (u-a_{U})|\right)\right)\\
	&=\eta(lf+h)e^{nu} + G,
  \end{split}
\end{equation}

Using the equation (\ref{soboelev bound 1}) and (\ref{soboelev bound 2}), we can show
$$||G||_{L^{q}(U)}\leq C_U \quad q\in \left[1,\frac{n}{n-1}\right).$$
Then by the H\"{o}lder inequality we can find $s_1>0$ such that
$$||G||_{L^{1}(B_{s_1}(x))}\leq \frac{1}{8}(n-1)!\gamma_n \quad x\in U.$$

Since $\int_M |he^{nu}|d\mu\leq\left(\int_Mh^\alpha e^{nu}d\mu\right)^{\frac{1}{\alpha}}\left(\int_Me^{nu}d\mu\right)^{\frac{\alpha}{\alpha-1}}$, we can find $A\in \mathbb{R}$ such that
$$\int_M |he^{nu}|d\mu\leq \frac{1}{8}(n-1)!,$$
for $E_\alpha\leq A$.
Then the second alternative in Proposition \ref{weak compactness} cannot hold.
As a consequence of Remark \ref{remark weak compactness}:
$$\int_U e^{nq_1w}d\mu\leq C_U$$
for some $q_1>1$. On the other hand, using the Jensen's formula, one can show $$e^{na_{U}}\leq C_U\int_U e^{nu}d\mu\leq C_U.$$
Since $w=u-a_{U}$ on $V$, we have  $$\int_V e^{nqu}d\mu\leq e^{nqa_{U}}\int_Ve^{nqw}d\mu\leq C_V.$$ Choosing $q_2=\frac{\alpha q_1}{\alpha+q_1-1}$, we can use the H\"older inequality to deduce  that:
$$\int_V (he^{nu})^{q_2}d\mu\leq\left(\int_V h^\alpha e^{nu}d\mu\right)^{\frac{q_1}{\alpha+q_1-1}}\left(\int_V e^{nq_1u}d\mu\right)^{\frac{\alpha-1}{\alpha+q_1-1}}\leq C.$$
Thus the right hand of the equation (\ref{equation in key proposition}) has a uniform $L^{q_2}(V)$ bound.
From standard elliptic estimates, we conclude
$$||w||_{W^{n,q_2}(W)}\leq C_W,$$
where $C_W$ is a uniform constant.  Using the Sobolev inequality we know $w$ is pointwise bounded in $W$. Then we can check that the right hand of the equation (\ref{equation in key proposition}) has a uniform $L^{\alpha}(W)$ bound.
Using elliptic estimates again , we can get
$$||u-a_{U}||_{W^{n,\alpha}(K)}=||w||_{W^{n,\alpha}(K)}\leq C $$
for some uniform constant $C$.
\hfill $\Box$ \\

\medskip

Observing the scaling invariance of $L^\infty$ norm, we have the following corollary:

\begin{corollary}\label{osc estimate}
Let $B_s(x)$ be a geodesic ball in M, and $u$ be the solution of the equation (\ref{sequence equation}). There is a uniform $A\in \mathbb{R}$ such that when $E_\alpha\leq A$, we have the following:
	
	If  $\int_{B_s(x)} lfe^{nu}d\mu\leq \frac{1}{8}(n-1)!\gamma_n $, then we have
	$$||u-a_{B_{ \frac{s}{2} } (x)}||_{L^\infty (B_{ \frac{s}{4}}(x)  )}\leq C,$$
	where $C=C(M, n, \alpha)$ is a constant and $a_{B_{ \frac{s}{2} } (x)}=\fint_{B_{ \frac{s}{2} } (x)}ud\mu$.

\end{corollary}
Inspired by Proposition \ref{strong compactness},  we give the following definition:
\begin{definition}\label{radius}
	Let $u$ is a solution of the equation (\ref{sequence equation}), we define the volume concentration radius function $s(x)$ to be the unique real function such that :
$$\int_{B_{s(x)}(x)}lfe^{nu}d\mu=\frac{1}{8}(n-1)!\gamma_n.$$
\end{definition}
Notice that $s(x)$ is a continous function of $x$ and it measure the scale of volume concentration.\\
In the case that $$\inf_{x\in M} s(x)\geq s $$
for some $s>0$, we can apply Proposition \ref{strong compactness} to the case  $U= K= M$, thus we get
$$||u-\bar{u}||_{W^{ n,\alpha } (M)}\leq C$$
for some constant $C= C(n, \alpha, s)$ . In view of the Sobolev embedding theorem, we know $|u-\bar{u}|\leq C$ pointwisely for some uniform constant $C$. Denoting $\int_Me^{nu}d\mu$ by $V$, we have
$$V=\int_Me^{n\bar{u}}e^{nu-n\bar{u}}d\mu \leq C\int_Me^{n\bar{u}}d\mu=Ce^{n\bar{u}}.$$
On the other hand, by the Jensen's formula,  we know
$$e^{n\bar{u}}\leq C\int_Me^{nu}d\mu= CV.$$
Thus $|\bar{u}|\leq C$ for some uniform $C$ . Combining with $||u-\bar{u}||_{W^{n,\alpha}(M)}\leq C$, we can conclude $||u||_{W^{n,\alpha}(M)}\leq C$ for some uniform $C$.

So when $s(x)$ has a uniform nozero lower bound, $u$ a has uniform $W^{n,\alpha}$ bound.

Now we care about the case where $\inf_{x\in M} s(x)$ is closed to $0$.

In this case , we need to describe the behavior of $u$ locally. And we shall rescale $u$ first.

For any given pair $(p,s)\in M \times \mathbb{R}^+$, We set
\begin{equation}\label{rescaled function}
	\hat{u}_{p,s}(x)=u(exp_{p}(sx))+\log s,\quad x\in B_{\delta/s}(0),
\end{equation}
where $B_{\delta/s}(0)\subset \mathbb{R}^n$ is the Euclidean ball of center $0$ and radius $\delta/s$, and  $\delta$ is a positive number less than the injectivity radius of $M$. We note that the ball $B_{\delta/s}(0)$ approaches
$\mathbb{R}^n$ when $s\to 0$. Let $T_{p,s}: B_{\delta/s}(0)\rightarrow M$ be defined by $T_{p,s}(x)=exp_{p}(sx)$, and define a metric on $B_{\delta/s}(0)$ by
$$g_{p,s} = s^{-2}T^*g.$$
Then $g_{p,s}\to g_0$ in $C^m(B_R(0))$ as $s\to 0$ for all $m\in\mathbb{N}$ and all $R>0$, where $g_0$ is the standard Euclidean metric of $\mathbb{R}^n$. An easy calculation shows that $\hat{u}_{p,s}$ satisfies the following equation:
\begin{equation}\label{rescaled discrete equation}
  P_{g_{p,s}}\hat{u}_{p,s}+s^n\hat{Q}_{p,s}=(l\hat{f}_{p,s}+ \hat{h}_{p,s})e^{n\hat{u}_{p,s}},
\end{equation}
where $P_{g_{p,s}}$ is the GJMS operator of the metric $g_{p,s}$ in $B_{\delta/s}(0)$ and
$$\hat{Q}_{p,s}(x)=Q(exp_{p}(sx)),\quad x\in B_{\delta/s}(0),$$
$$\hat{f}_{p,s}(x)=f(exp_{p}(sx)),\quad x\in B_{\delta/s}(0),$$
$$\hat{h}_{p,s}(x)=h(exp_{p}(sx)),\quad x\in B_{\delta/s}(0),$$
and we also have $\int_{B_{\delta/t}(0)}\hat{h}^\alpha_{p,s}e^{n\hat{u}_{p,s}}d\mu_{g_{p,s}}=\int_{B_{\delta}(p)}h^\alpha e^{nu}d\mu$ from a direct calculation, where $d\mu_{g_{p,s}}$ is the volume element associated to  $g_{p,s}$.

Noticing that $\int_{B_{\delta/s}(0)}e^{n\hat{u}_{p,s}}d\mu_{g_{p,s}} =\int_{B_{\delta}(p)}e^{nu}d\mu$,  we have $$\int_{B_1(0)}l\hat{f}_{p,s}e^{n\hat{u}_{p,s}}d\mu_{g_{p,s}}=\frac{1}{8}(n-1)!\gamma_n.$$
If we choose $p$ as one of the global minimum point of $s(x)$, and $s=s(x)$, we can also have

\begin{equation}\label{volume bdd}
	\int_{B_1(x)}l\hat{f}_{p,s}e^{n\hat{u}_{p,s}}d\mu_{g_{p,s}}\leq \frac{1}{8}(n-1)!\gamma_n,\quad x\in B_{\delta/s}(0)
\end{equation}
We have the following compactness result on $\hat{u}_{p,s}$:
\begin{proposition}\label{bubble convergence}
	Assume $u$ is a solution of the equation (\ref{sequence equation}),  $p$ is one of the global minimum point of $s(x)$, and $s=s(p)=\inf_{x\in M}s(x)$.
	Then for $E_\alpha$ and $s$ are both small enough (that is, smaller than a given small real number), there exist $\varepsilon>0$ and  $x_0\in \mathbb{R}^n$ such that,
	\begin{equation}\label{bubble app}
	 \|\hat{u}_{p,s} - U_{x_0, \varepsilon} \|_{ W^{ n,\alpha } (B_{ L(E_p, s) } (0)) }\leq J(E_\alpha, s)  ,
	\end{equation}
  where the function $\hat{u}_{\infty} $ is given by
  \begin{equation}\label{bubble}
 	U_{x_0, \varepsilon}=-\log(1+\frac{\lambda^2}{\varepsilon^2}|x-x_0|^2)-\log\varepsilon, \qquad \lambda=\left(\frac{lf(p)}{2^n(n-1)!}\right)^{\frac{1}{n}},
 \end{equation}
$J$ and $L$ are positive numbers satisfy
\[
	\lim_{ E_\alpha, s\to 0 } J(E_\alpha, s)=0 \qquad \lim_{ E_\alpha, s\to 0 } L(E_\alpha, s)=\infty
.\] Moreover, there exists a positive real number $J_1(M, E_\alpha, s)$ , such that:
  \begin{equation}\label{volume of bubble}
  \left| \int_{B_{Ms}( p ) }lfe^{nu}d\mu - (n-1)!\gamma_n \right| \leq J_1(M, E_\alpha, s) ,
  \end{equation}
  where $J_1(M, E_\alpha, s)$ is a positive real number satisfies
  \[
	  \lim_{M\to\infty}\lim_{E_\alpha, s\to 0 }J_1(M, E_\alpha, s)= 0
  .\]
  We also have $\varepsilon\leq c(n)$ for some uniform constant $c(n)>0$, when $E_\alpha$ and $s$ are both small enough.
\end{proposition}
\medskip

{\bf Proof. } We argue by contradiction. Assume there exists a sequence of function $\{u_k\}_{k\in \mathbb{N}}$ satisfy:
\begin{equation}
  \begin{cases}
  Pu_k +Q = (l_kf_k+ h_k)e^{nu_k},
\\
\int_M l_kf_ke^{nu_k}d\mu = \int_M Q d\mu = (n-1)!\gamma_n
\\
\lim_{k\to \infty} E_{\alpha,k} =\lim_{k\to\infty} \int_M g_k^\alpha e^{nu_k}d\mu = 0
,
  \end{cases}
\end{equation}
and we denote the volume concentration radius function of $u_k$ by $s_k(x)$. Then we can choose $p_k$ as a minimum point of $s_k(x)$ and let $s_k= s_k(x_k)$. We write $\hat{u}_k=\hat{u}_{p_k,s_k}$ as the rescaled function of $u_k$.

We now assume there exist a $L_1>0$ and $J_1>0$ such that
\[
\|\hat{u}_k - U_{x_0, \varepsilon} \|_{ W^{ n,\alpha } (B_{ L_1 } (0)) }\geq J_1
,\]
for any $k$, $x_0$ and $\varepsilon$, since otherwise (\ref{bubble app}) will hold for all solution $u$ of the equation (\ref{sequence equation}).

By Proposition \ref{weak global estimate} and some scaling argument one finds $$\int_{ B_R(0) }
|\nabla^j_{ g_k } \hat{u}_k|^qd\mu\leq C_U\quad q\in\left[1,n/j\right),\quad j\in [1, n-1]\cap\mathbb{N}$$
for any $R>0$. Then arguing as in Proposition \ref{strong compactness}  we obtain $||\hat{u}_k-{\hat{a}_{k,B_R(0)}}||_{W^{n,\alpha}\left( B_R(0) \right) }\leq C_R$, where $C_R$ is independent of $k$ and $\hat{a}_{k,B_R(0)}=\fint_{B_{R}(0)}\hat{u}_kd\mu_{g_k}$.

From the Jensen's inequality, we can deduce
$$e^{n{\hat{a}_{k,B_R(0)}}}\leq C_R\int_{B_R(0)}e^{n\hat{u}_k}d\mu_{g_k}\leq C_R.$$
By a direct calculation, we have that
$$\frac{1}{8}(n-1)!\gamma_n\leq \int_{B_R(0)}e^{n{\hat{a}_{k,B_R(0)}}}e^{nu-{\hat{a}_{k,B_R(0)}}}d\mu_{g_k} \leq C_R\int_{B_R(0)}e^{n{\hat{a}_{k,B_R(0)}}}d\mu_{g_k}=C_Re^{n{\hat{a}_{k,B_R(0)}}}.$$
Then we can conclude $||\hat{u}_k||_{W^{ n, \alpha } \left(  {B_R(0)}\right) }\leq C_R$ for some $C_R$ independent of $k$. So after passing to a subsequence, we can assume $\hat{u}_k\to\hat{u}_{\infty}$ in the local weak  $W^{n,\alpha}$ topology, where $\hat{u}_{\infty}$ satisfies the equation
\begin{equation}\label{bubble equation}
  \begin{cases}
	  \left(  -\Delta_0\right) ^{ \frac{n}{2} }  \hat{u}_{\infty}= l_{\infty}f(p)e^{n\hat{u}_{\infty}}\\
    \int_{\mathbb{R}^n}e^{n\hat{u}_{\infty}}dx\leq V.
  \end{cases}
\end{equation}
Moreover, using the local $W^{n,\alpha}$ bound we can rewrite the equation (\ref{rescaled discrete equation}) as the following:
$$\left(  -\Delta_0\right) ^{ \frac{n}{2} } \hat{u}_k= l_{\infty}f(p)e^{n\hat{u}_k}+ H_k,$$
where $H_k$ satisfying $\lim_{k\to+\infty}\int_{B_R(0)}H_k^\alpha\mu_0=0$ for any $R>0$. As a consequence, we can improve the
local weak convergence of $\hat{u}_k$ to the local strong convergence by invoking elliptic estimates.

The solutions of the equation (\ref{bubble equation}) have been classified in \cite{WX1999},\cite{L2001} and \cite{M20092} which tells us $\hat{u}_{\infty}$ has the form
$$\hat{u}_{\infty}=-\log(1+\frac{\lambda_\infty^2}{\varepsilon^2}|x-x_0|^2)-\log\varepsilon, \qquad \lambda_\infty=\left(\frac{l_{\infty}f(p)}{2^n(n-1)!}\right)^{\frac{1}{n}}.$$
That is , $\hat{u}_\infty= U_{x_0,\varepsilon}$ for some $x_0$ and $\varepsilon$.
Since $\hat{u}_k$ converges to $\hat{u}_\infty$ strongly in $W^{n,\alpha}_{local}(\mathbb{R}^n)$. And this is contract to our assumption. So we finish the proof of (\ref{bubble app}).

By a simple calculation we have $\int_{\mathbb{R}^n}l_{\infty}f(p)e^{n\hat{u}_{\infty}}dx= (n-1)!\gamma_n$. As a consequence, the inequality (\ref{volume of bubble}) follows from (\ref{bubble app}).
The last assertion of this proposition comes from the inequality (\ref{volume bdd})  .
\hfill $\Box$ \\

\medskip

\begin{corollary}\label{stdbubble convergence}
	Assume $u$ is a solution of the equation (\ref{sequence equation}),
	Then for $E_\alpha$ and $s_{min}=\inf_{x\in M}s(x)$ are both small enough (that is, smaller than a given small real number), there exist $s>0$ and  $p\in M$ such that,
	\begin{equation}\label{stdbubble app}
	 \|\hat{u}_{p,s} - U_{0, 1} \|_{ W^{ n,\alpha } (B_{ L(E_\alpha, s_{min}) } (0)) }\leq J(E_\alpha, s_{min})  ,
	\end{equation}
  where the function $\hat{u}_{\infty} $ is given by
  \begin{equation}\label{stdbubble}
 	U_{0, 1}=-\log(1+\lambda^2|x-x_0|^2), \qquad \lambda=\left(\frac{l_{\infty}f(p)}{2^n(n-1)!}\right)^{\frac{1}{n}},
 \end{equation}
and $J$, $L$ are positive numbers satisfy
\[
	\lim_{ E_\alpha, s_{min}\to 0 } J(E_\alpha, s_{min})=0 \qquad \lim_{ E_\alpha, s_{min}\to 0 } L(E_\alpha, s_{min})=\infty
.\] Moreover, there exists a positive real number $J_1(M, E_\alpha, s_{min})$ , such that:
  \begin{equation}\label{volume of stdbubble}
	  \left| \int_{B_{Ms}( p ) }lfe^{nu}d\mu - (n-1)!\gamma_n \right| \leq J_1(M, E_\alpha, s_{min}) ,
  \end{equation}
  where $J_1(M, E_\alpha, s_{min})$ is a positive real number satisfies
  \[
	  \lim_{M\to\infty}\lim_{E_\alpha, s_{min}\to 0 }J_1(M, E_\alpha, s_{min})= 0
  .\]
  We also have $s_{min}\leq s\leq c(n)s_{min}$.
 	
\end{corollary}

\medskip
{\bf Proof. }
We denote the $p$ and $s$ in Proposition \ref{bubble convergence} by $p_{min}$ and $s_{min}$, and choose
\[
	p=exp_{p_{min}}(x_0)\qquad\qquad s=\varepsilon s_{min}
.\]
Then this corollary follows easily.

\hfill $\Box$ \\

\medskip

From now we fix $p$ and $s$ such that they fulfill the conclusion of Corollary \ref{stdbubble convergence}.

We now describe the behavior of $u$ in $M\setminus \{p\}$ when  $E_\alpha$ and $s$ are both tend to $0$ .

\begin{proposition}\label{outside convergence}
 Assume $u$ is a solution of the equation (\ref{sequence equation}) and $E_\alpha$, $s$  as above.
	Then for $E_\alpha$ and $s$ are both small enough (that is, smaller than a given small real number), we have:
	$$\|u-\bar{u}- G_p\|_{W^{n,\alpha} \left(M\setminus B_{ \delta(E_\alpha,s) } (p)\right) }  \leq J(E_\alpha, s),$$
  where $G_p$ is the Green's function satisfying
  $$\begin{cases}
	  PG_p+Q=(n-1)!\gamma_n\delta_p\\
    \int_M G_pd\mu=0.
  \end{cases}$$
  in the distribution sense. And $\delta(E_\alpha,s)$ and $J(E_\alpha, s)$ are positive real numbers such that
  \[
	  \lim_{ E_\alpha,s\to 0 } \delta(E_\alpha,s)=0\qquad \lim_{ E_\alpha,s\to 0 } J(E_\alpha, s)=0
  .\]

\end{proposition}
\medskip
{\bf Proof. }
We argue by contradiction. So we can assume there exists a sequence of function $\{u_k\}_{k\in \mathbb{N}}$ satisfy:
\begin{equation}
  \begin{cases}
  Pu_k +Q = (l_kf_k+ h_k)e^{nu_k},
\\
\int_M l_kf_ke^{nu_k}d\mu = \int_M Qd\mu = (n-1)!\gamma_n
\\
\lim_{k\to \infty} E_{\alpha,k} =\lim_{k\to\infty} \int_M g_k^\alpha e^{nu_k}d\mu = 0
.
  \end{cases}
\end{equation}
and there exists $\delta_1>0$ and $J_1>0$ such that
\[
	\|u_k-\bar{u}_k- G_{p_k}\|_{W^{n,\alpha} \left(M\setminus B_{ \delta_1 } (p_k)\right) }  \geq J_1
.\]
By Proposition \ref{weak global estimate}, using the compactness of $M$ and the Poincar\'{e} inequality , we have
$$||u_k-\bar{u}_k||_{W^{n-1,q}(M)}\leq C_q,\quad q\in\left(1,\frac{n}{n-1}\right)$$
for some $C_q$ independent of $k$.

Then after passing to a subsequence, we can assume $u_k-\bar{u}_k\to u_\infty$ weakly in $W^{n-1,q}(M)$ for some $u_\infty\in W^{n-1,q}(M)$.

In view of the formula (\ref{volume of bubble}), we know that
$$\int_{K}l_kf_ke^{nu_k}d\mu\to 0$$
for any compact set $K\subset M\setminus\{p\}$ as $k\to+\infty$. Then we can apply Proposition \ref{strong compactness}
to get the following estimate:
$$||u_k-a_{k,K}||_{W^{n,\alpha}(K)}\leq C_K,$$
where $C_K$ is independent of $k$. Thus we obtain
$$|a_{k,K}-\bar{u}_k|\leq C_K\left(||u_k-a_{k,K}||_{W^{n,\alpha}(K)}+||u_k-\bar{u}_k||_{W^{n-1,q}(M)} \right)\leq C_K.$$
So we can conclude that $||u_k-\bar{u}_k||_{W^{n,\alpha}(K)}\leq C_K$. By passing to a subsequence we can deduce
$$u_k-\bar{u}_k\to u_\infty\ \qquad in\quad W_{loc}^{n,\alpha}(M\setminus\{p\} )$$
On the other hand, from the equation (\ref{sequence equation}), we know $u_k-\bar{u}_k$ satisfies
$$P(u_k-\bar{u}_k) +Q =e^{n\bar{u}_k}(l_kf_k+h_k)e^{nu_k-n\bar{u}_k}.$$
Let $k\to+\infty$ we get
$$\begin{cases}
	Pu_\infty+Q=(n-1)!\gamma_n\delta_p\\
  \int_M u_\infty d\mu=0.
\end{cases}$$

So we can conclude $u_\infty= G_p$. This is contract to our assumption and we finish the proof of this proposition.
\hfill $\Box$ \\

In view of Proposition \ref{stdbubble convergence} and Proposition \ref{outside convergence}, we know the behavior of $u$ on $B_{Ls}(p)$ and $M\setminus B_{\delta}(p)$ pretty well.

To understand the behavior on $B_\delta(p)\setminus B_{Ls}(p)$ , we need the following two lemmas:

\begin{lemma}\label{one bubble}
	Assume $u$ is a solution of the equation (\ref{sequence equation}). Then for $E_\alpha$ small enough, there exists a uniform constant $C>0$, such that
	\[
		d_g(x,p)\leq C s(x)
	.\]
\end{lemma}
{\bf Proof. }
By choose $M$ large enough and $\beta>0$ small enough, we can deduce from the formula (\ref{volume of stdbubble})

\[
	\int_{B_{Ms}(p)}lfe^{nu}d\mu\geq \frac{7}{8}(n-1)!\gamma_n
,\]
for $s\leq \beta$ and $E_\alpha\leq \beta$.
Then the definition of $s(x)$ tells us
\[
	s(x)\geq d_g(x,p)-Ms
.\]
So when $d_g(x,p)\geq 2Ms$, we have
\[
	s(x)\geq Ms\geq\frac{d_g(x,p)}{2}
.\]
When $d_g(x,p)\leq 2Ms$, we know from $s(x)\geq s_{min}\geq\frac{s}{c(n)}$ that
\[
	s(x)\geq\frac{d_g(x,p)}{2Mc(n)}
	.\]
We now can set $C=max\{2, 2Mc(n), Diam(M)/\beta\}$

\hfill $\Box$ \\

\begin{lemma}\label{averdecay on neck}
	Assume $u$ is a solution of the equation (\ref{sequence equation}) and $\nu\in [1, 2)$. Then for $E_\alpha$ and $s$ are both small enough ( that is, $E_\alpha\leq \beta$, $s\leq \beta$ for some small real number $\beta>0$ to be determined), there exist $\delta>0$, $L>0$ independent of $s$, such that
	\(
		r^{n\nu}e^{n\bar{u}_p(r)}
	\)
	is monotonically decreasing on $[Ls, \delta]$, where
	\[
		\bar{u}_p(r)= Vol\left( \partial B_r(p) \right) ^{-1}\int_{ \partial B_r(p) } u(x) d\sigma
	.\]
\end{lemma}
{\bf Proof. }
We argue by contradiction.  Assume there exists a sequence of function $\{u_k\}_{k\in \mathbb{N}}$ which satisfy:
\begin{equation}
  \begin{cases}
  Pu_k +Q = (l_kf_k+ h_k)e^{nu_k},
\\
\int_M l_kf_ke^{nu_k}d\mu = \int_M Q d\mu = (n-1)!\gamma_n
\\
\lim_{k\to \infty} E_{\alpha,k} =\lim_{k\to\infty} \int_M g_k^\alpha e^{nu_k}d\mu = 0
.
  \end{cases}
\end{equation}
And $p_k$, $s_k$ is the parameter of $u_k$ determined by Proposition \ref{stdbubble convergence}.
\[
	A_k(r)=\nu r+\bar{u}_{p_k}(r)
.\]
We know from Proposition \ref{stdbubble convergence} that
$s_kA_k^\prime(s_kr)$ converges to $\frac{\nu}{r}-\frac{2\lambda^2r}{1+\lambda^2r^2}$ in $W_{loc}^{n-1,\alpha}(\mathbb{R}^+)$, thus in $C_{loc}^0(\mathbb{R}^+)$. And this implies that for any $R\geq 2R_\nu : =\frac{2}{\lambda}\sqrt{\frac{\nu}{2-\nu}}$, there exists $k_0(R)$ such that
\[
	A_k^\prime(s_kr)<0 \text{ for } k\geq k_0(R), \quad r\in [2R_\nu, R]
.\]
Define
\[
	r_k=sup \left\{ r\in [2s_kR_\nu, \delta_0]: A^\prime(r)<0 \text{ in } [2s_kR_\nu, r) \right\}
.\]
From above discussion we infer that
\[
	\lim_{k\to\infty}\frac{r_k}{s_k}=\infty
.\]
 Now we can assume that
\[
	\lim_{k\to\infty}r_k=0
.\]
Otherwise, $\{u_k\}_{k\in \mathbb{N}}$ will obey the consequence of this Lemma.

Consider
\[
	v_k(x)= \widehat{\left( u_k \right)}_{p_k,s_k}- C_k
,\]
where $C_k$ is a constant such that
\[
	\int_{\partial B_1(0)}v_k d\sigma_k =0
.\]
Then we have
\begin{equation}\label{critical radius}
	\frac{d}{dr}\left( r^{n\nu}e^{n\bar{v}_k} \right) (1)=0.
\end{equation}

Argue as in Proposition \ref{outside convergence}
, after passing to a subsequence,
\[
	v_k\longrightarrow G_0 \quad\text{ in}\quad W_{loc}^{n, \alpha}\left( \mathbb{R}^n\setminus\{0\}\right)\quad \text{as}\quad k\to\infty
.\]
Where $G_0$ satisfies
\[
	\left( -\Delta_0 \right) ^{\frac{n}{2}}G_0 = (n-1)!\gamma_n
.\]
Moreover, thanks to Corollary \ref{osc estimate} and Lemma \ref{one bubble} , we have
\[
	\max_{B_\frac{|x|}{4C}(x)}G_0-\min_{B_\frac{|x|}{4C}(x)}G_0\leq C(M,n,\alpha)
.\]
Then we know from standard elliptic estimates:
\[
	\left| \nabla_0^k G_0 \right| (x)\leq \frac{C(M,n,k,\alpha)}{|x|^k}
,\]
hold for all $k\in \mathbb{N}$.
So
\[
G_0 = 2\log\frac{1}{|x|} + C
.\]
Since $\nu<2$, we get in particular that
\[
	\frac{d}{dr}\left( r^{n\nu}e^{n\bar{G}_0} \right) (1)<0
.\]
And this is contract to (\ref{critical radius}).

\hfill $\Box$ \\

\begin{proposition}\label{decay on neck}
	Assume $u$ is a solution of the equation (\ref{sequence equation}). Then for $E_\alpha$ and $s$ are both small enough, there exist $\delta>0$, $L>0$ and $\eta>0$ such that
	\[
		d_g(x, p)^{n\nu}e^{nu(x)}\leq \eta s^{n\left( \nu-1 \right) } \qquad \text{for } x\in B_\delta(p)\setminus B_{Ls}(p)	
	.\]

\end{proposition}
{\bf Proof. }
This result follows immediately from combing Proposition \ref{averdecay on neck}, Corollary \ref{osc estimate} and Lemma \ref{one bubble}.  \hfill $\Box$ \\

Proposition \ref{stdbubble convergence} and Proposition \ref{decay on neck} allow us to get control of $\hat{u}_{p,s}$ on the whole $B_{\frac{\delta}{s}}(0)$.

\begin{corollary}\label{global upp}
		Assume $u$ is a solution of the equation (\ref{sequence equation}). Then for $E_\alpha$ and $s$ are both small enough, there exists a constant $C\geq 0$ independent of $s$ such that
		\[
			e^{n\hat{u}_{p,s}}\leq C \text{ on } B_{\frac{\delta}{s}}(0)
		.\]

In particular, we have
\[
	\int_{B_{\frac{\delta}{s}}(0)}\hat{h}_{p,s}^\alpha e^{\alpha n\hat{u}_{p,s}}d\mu_{g_{p,s}}\leq C\int_{B_{\frac{\delta}{s}}(0)}\hat{h}_{p,s}^\alpha e^{n\hat{u}_{p,s}}d\mu_{g_{p,s}}
.\]
\end{corollary}
Combing Proposition \ref{decay on neck} and Proposition \ref{outside convergence}, we have the following corollary.
\begin{corollary}\label{out vol}
	Assume $u$ is a solution of the equation (\ref{sequence equation}) and $\nu\in [1, 2)$. Then for $E_\alpha$ and $s$ are both small enough, there exists a constant $C>0$ independent of $s$, such that
	\[
		\bar{u}\leq \left( \nu-1 \right) \log s + C
	.\]
In particular, we have
\[
	\int_{M\setminus B_\delta(p)}e^{nu}d\mu\leq C_1 s^{n\left( \nu-1 \right)},
\]
for some constant $C_1$ independent of $s$.
\end{corollary}

\section{The Test Function $\varphi_{p,s}$}

\setcounter{equation}{0}

\medskip

In this section, we will construct the test function $\varphi_{p,s}$ which plays a important role in the proof of our results. Our construction extends the test function in \cite{DJLW1997} and \cite{LLL2012} to the general even dimension. A similar test function has been used by C.B. Ndiaye  in \cite{N2019}.

Using the existence of conformal normal coordinates (see \cite{C1991}, \cite{G1993}), we have that for any $p\in M$ there exists a function $u_p\in C^\infty(M)$ such that

\[
	g_p= e^{2u_p}g \quad \text{ verifies } \quad detg_p(x)=1 \quad \text{ for } \quad x\in B
,\]
with $0<\rho_p<\frac{inj_{g_p}(M)}{10}$. Moreover, we can take the families of functions $u_p$, $g_p$ and $\rho_p$ such that
\[
	\text{the map } p\longrightarrow u_p \quad\text{ are }\quad C^1 \quad \text { and } \rho_p\geq\rho_0>0
,\]
for some small positive $\rho_0$ satisfying $\rho_0<\frac{inj_g}{10}$, and
\[
	\|u_p\|_{C^n(M)}=\mathcal{O}(1),\quad \frac{1}{C^2}g\leq g_p\leq C^2 g
,\]
\[
	u_p(x)=\mathcal{O}\left(d_{g_p}^2(p,x)\right)=\mathcal{O}\left( d_{g_p}^2(p,x) \right)
,\]
\[
	u_p(p)=0,\qquad \nabla u_p (p) = 0, \qquad R_{g_p}(p)=0,\qquad Ric_{g_p}(p)=0
,\]
for some large positive constant $C$ independent of $p$. Furthermore, using $\nabla u_p (p)= 0$ and the scalar curvature equation, namely
\[
	-\frac{4(n-1)}{n-2}\Delta_{g_p}\left( e^{-\frac{n-2}{2}u_p} \right) + R_{g_p}e^{-\frac{n-2}{2}u_p} = R_{g}e^{-\frac{n+2}{2}u_p}
,\]
for $n\geq 3$,
and
\[
	2\Delta_{g_p}(u_p) + R_{g_p} = R_g e^{-2u_p}
,\]
for $n=2$,
it is easy to see that the following holds

\begin{equation}\label{local delta}
	\Delta_{g_p}u_p(p)=\frac{R_{g}(p)}{2(n-1)}
\end{equation}
On the other hand, using the properties of $g_p$, it is easy to check that for every $u\in C^2(M)$ there holds
\[
	\nabla_{g_p}u(p)=\nabla u(p)=\nabla_{\mathbb{R}^n}\hat{u}(0),\qquad \Delta_{g_p}u(p)=\Delta u(p)= \Delta_{\mathbb{R}^n}\hat{u}(0)
,\]
where
\[
	\hat{u}(y)= u\left( exp_p^{g_p}(y) \right) ,\quad y\in \mathbb{R}^n
.\]

We have
\begin{equation}
  \begin{split}
	  \Delta_{g_p} f(r_p)
&=f^{\prime\prime}(r_p)+(n-1)\frac{f^{\prime}(r_p)}{r_p},\\
  \end{split}
\end{equation}
where $r_p(x)=d_{g_p}(x,p)$

Using the conformal normal coordinate, we can give a simple proof of the following expansion of the Green's function.

\begin{proposition}\label{Green expan}
	  Assume $G_p$ is the solution of the equation (\ref{Green}). Then
   \[
	   G_p(x)=H_p(x)+K_p(x)
   \]
   is smooth on $M\setminus \{p\}$, $H_p$ extends to a twice differentiable function on $M$ and
   \[
	   K_p(x)=-2f(r)\log r
   ,\]
   where $r(x)=d_g(p,x)$ is the geodesic distance from $x$ to $p$, $f(r)$ is a $C^\infty$ positive function such that $f(r)=1$ on a neighborhood of $r=0$, and $f(r)=0$ for $r\geq inj_g(M)$.
  \end{proposition}

{\bf Proof. }
We denote $r_p(x)=d_{g_p}(p,x)$ . Then the knowledge from Riemannian geometry tells us:
\begin{equation}\label{expansion of dist}
	r_p^2(x)= r^2(x) + \frac{1}{3}\left(  \nabla^2u_p(p)\right) (X,X)r^2(x)+\mathcal{O}(r^5),
\end{equation}
where $x=exp_p(X)$. So we know that $\frac{r_p^2}{r^2}$ is twice differentiable at $p$.

Since the GJMS operator is a natural operator in Riemannian geometry, it will have the following form:
\[
P= \left( -\Delta \right) ^{ \frac{n}{2} } + aR\left( -\Delta \right) ^{ \frac{n-2}{2} } + bRic_{ij}\nabla_i\nabla_j\left( -\Delta \right) ^{ \frac{n-4}{2} } + l.o.t
.\]
Here $R$ is the scalar curvature, and $Ric$ is the Ricci curvature tensor.

Since $R_{g_p}(p)=0$ and $Ric_{g_p}(p)=0$, we have
\[
	P_{g_p}\left(-2\log r_p \right)  = \mathcal{O}\left( \frac{1}{r_p^{n-3}} \right) \in L^q \text{ for any } q<\frac{n}{n-3}
.\]
On the other hand, from the conformally covariance, we have
\[
	P_{g_p}\left( G_p-u_p \right) +Q_{g_p}=(n-1)!\gamma_n\delta_p
.\]
The by basic elliptic estimate, we deduce
\[
	G_p-u_p+2\log r_p\in W^{n, q}\subset C^{2, n-2-\frac{n}{q}} \text{ for any } q<\frac{n}{n-3}
.\]
That is
\[
G_p + 2\log r + \left( \log\frac{r_p^2}{r^2}-u_p \right) \in C^{2, n-2-\frac{n}{q}} \text{ for any } q<\frac{n}{n-3}
.\]
Since $\log\frac{r_p^2}{r^2}-u_p$ is twice differentiable at $p$, the proof is finished.
 \hfill $\Box$ \\

 The expansion (\ref{expansion of dist}) also tells us that:
 \begin{equation}\label{coo diff}
	\Delta \left( \log\frac{r_p^2}{r^2}-u_p \right) = \frac{1}{3}\Delta u_p= \frac{R_g(p)}{6(n-1)}
 \end{equation}

We define the function as

\[
	\varphi_{p,s}=G_p+\eta_{p,s}\left( -\log\left( 1+\frac{s^2}{\lambda^2r_p^2} \right)+\sum_{k=0}^{\frac{n-2}{2}} B_{k,s}(r_p^2-L^2s^2)^{k}  \right)
.\]
Here $\eta_{p,s}$ is a smooth function such that $\eta_{p,s}=1$ on $B_{\frac{Ls}{3}}^{g_p}(p)$ and $\eta_{p,s}=0$ on $M\setminus B_{\frac{Ls}{2}}^{g_p}(p)$ with $|\nabla^k\eta_{p,s}|\leq\frac{C(n,k)}{L^ks^k}$,
and
\[
	B_{k,s}=\frac{1}{k!}\frac{d^k}{d r}\log\left( 1+\frac{s^2}{\lambda^2r^2} \right)\Bigg|_{r=Ls}
.\]
By a direct calculation, we find
\[
	\left| B_{k,s} \right|\leq\frac{C(k)}{L^{k+2}s^k}
.\]
Then we have
\[
	\left|\nabla^k\left( \eta_{p,s}(\phi_{p,s}-G_p)\right) \right|\leq  \frac{C(n,k)}{L^{k+2}s^k}
.\]
So we can calculate that
\begin{equation}\label{app phi}
\begin{split}
	&\left| P\varphi_{p,s}+ Q -\frac{lf(p)s^n\eta_{p,s}e^{nu_p}}{(s^2+\lambda^2r_p^2)^n} \right|\\
	=&e^{nu_p}\left| P_{g_p}\varphi_{p,s}+ Q_{g_p} -\frac{lf(p)s^n\eta_{p,s}}{(s^2+\lambda^2r_p^2)^n} \right|\\
	\leq &\eta_{p,s}e^{nu_p}  \left|\left( P_{g_p}- (-\Delta_{g_p})^{\frac{n}{2}} \right)\log \left( 1+\frac{s^2}{\lambda^2r_p^2 }\right)    \right|+ \eta_{p,s}\left| \sum_{k=0}^{\frac{n-2}{2}} B_{k,s}P_{g_p}(r_p^2-L^2s^2)^{k}  \right| +\frac{C\chi_{B_{\frac{Ls}{2}}(p)\setminus B_{\frac{Ls}{3}}(p)}}{L^{n+2}s^n}\\
	\leq &C \eta_{p,s}\left( r_p\left| \nabla_{g_p}^{n-2}\log \left( 1+ \frac{s^2}{\lambda^2r_p^2} \right) \right|+\sum_{k=1}^{n-3}\left| \nabla_{g_p}^k\log \left( 1+\frac{s^2}{\lambda^2r_p^2} \right) \right|    +\frac{1}{L^{n-1}s^{n-3}} \right)+\frac{C\chi_{B_{\frac{Ls}{2}}(p)\setminus B_{\frac{Ls}{3}}(p)}}{L^{n+2}s^n}\\
\leq & \frac{Cs^2\eta_{p,s}}{r_p^{n-1}}+\frac{C\chi_{B_{\frac{Ls}{2}}(p)\setminus B_{\frac{Ls}{3}}(p)}}{L^{n+2}s^n}.	\end{split}
\end{equation}

We now set $\mathcal{E}(x)=\sum_{k=0}^{\frac{n-2}{2}} B_{k,s}(r_p^2-L^2s^2)^{k}(x)$, and compute $E[\varphi_{p,s}]$.
We have
\[
\begin{split}
	&\frac{n}{2}\int_M \varphi_{p,s}P\varphi_{p,s}d\mu+ n\int_M Q\varphi_{p,s}d\mu\\
	=&\frac{n}{2}\int_M\varphi_{p,s}(PG_p+Q)d\mu+\frac{n}{2}\int_M(\varphi_{p,s}-G_p)(P\varphi_{p,s}+Q)d\mu+\frac{n}{2}\int_MQG_pd\mu\\
	=&\frac{n!\gamma_n}{2}\left( \log\frac{\lambda^2}{s^2} + H_p(p) + \mathcal{E}(p) \right)+\frac{n}{2}\int_MQG_pd\mu
	-\frac{n}{2}\int_{B_{Ls}(p)}\log\left( 1+\frac{s^2}{\lambda^2r_p^2}\right) \frac{lf(p)s^n}{(s^2+\lambda^2r_p^2)^n}d\mu_{g_p}\\
	+&\frac{n}{2}\int_{B_{Ls}(p)}\frac{lf(p)s^n\mathcal{E}}{(s^2+\lambda^2r_p^2)^n}d\mu_{g_p}+\mathcal{O}\left(s^3L\log\frac{1}{s} +\frac{1}{L^{n+2}}+\frac{1}{L^4} \right)\\
	=&\frac{n!\gamma_n}{2}\left( \log\frac{\lambda^2}{s^2} + H_p(p) + 2\mathcal{E}(p) -\frac{C_{\log}}{C_0}\right)+\frac{n}{2}\int_MQG_pd\mu+\mathcal{O}\left(s^3L\log\frac{1}{s} + \frac{1}{L^{n+2}}+\frac{1}{L^4} \right),
\end{split}
\]
where
\[
	C_{\log}=\int_{\mathbb{R}^n}\log\left( 1+\frac{1}{\lambda^2|x|^2} \right) \frac{1}{(1+\lambda^2|x|^2)^n}dx
,\]
and
\[
	C_0=\int_{\mathbb{R}^n} \frac{1}{(1+\lambda^2|x|^2)^n}dx
,\]
Notice that
\[
	\frac{C_{\log}}{C_0}=\int_{\mathbb{R}^n}\log\left( 1+\frac{1}{|x|^2} \right) \frac{1}{(1+|x|^2)^n}dx\Bigg/\int_{\mathbb{R}^n} \frac{1}{(1+|x|^2)^n}dx
\]
are constants depending only on $n$.
We can also calculate:
\[
\begin{split}
         &\lambda^{-2n}s^n\int_M fe^{n\varphi_{p,s}}d\mu\\
	=&\lambda^{-2n}s^n\int_{B_{Ls}(p)} fe^{n\varphi_{p,s}}d\mu + \mathcal{O}\left( \frac{1}{L^n} \right) \\
	=& \int_{B_{Ls}(p)}\frac{s^n}{(s^2+\lambda^2r_p^2)^n}e^{nH_p+n\mathcal{E}+\log f +n\log\frac{r_p^2}{r^2}-nu_p}d\mu_{g_p}+\mathcal{O}\left( \frac{1}{L^n} \right) \\
	=&\int_{B_{Ls}(p)}\frac{s^n}{(s^2+\lambda^2r_p^2)^n}e^{nH_p+n\mathcal{E}+\log f +n\log\frac{r_p^2}{r^2}-nu_p}(p)d\mu_{g_p}\\
     +&\frac{1}{n}\int_{B_{Ls}(p)}\frac{s^nr^2}{(s^2+\lambda^2r_p^2)^n}\Delta_{g_p}e^{nH_p+n\mathcal{E}+\log f +n\log\frac{r_p^2}{r^2}-nu_p}(p)d\mu_{g_p}\\
	+&\mathcal{O}\left(s^3+\frac{1}{L^5}+\frac{1}{L^n}   \right)\\
	=&\int_{B_{Ls}(p)}\frac{s^n}{(s^2+\lambda^2r_p^2)^n}e^{nH_p+n\mathcal{E}+\log f +n\log\frac{r_p^2}{r^2}-nu_p}(p)d\mu_{g_p}\\
+&\frac{1}{n}\int_{B_{Ls}(p)}\frac{s^nr^2}{(s^2+\lambda^2r_p^2)^n}\Delta_{g_p}e^{nH_p+\log f +n\log\frac{r_p^2}{r^2}-nu_p}(p)d\mu_{g_p}\\
	+&\mathcal{O}\left(s^3+\frac{1}{L^4} +\frac{1}{L^n}  \right).\\
\end{split}
\]
We denote
\[
	C_2=\int_{\mathbb{R}^n}\frac{|x|^2}{\left(  1+\lambda^2|x|^2\right)^n }dx
.\]
Then for $n>2$, we have
\[
	\begin{split}
	&\lambda^{-2n}s^n\int_{B_{Ls}(p)} fe^{n\varphi_{p,s}}d\mu\\
		=&C_0 e^{nH_p+n\mathcal{E}+\log f +n\log\frac{r_p^2}{r^2}-nu_p}(p)+\frac{C_2s^2}{n}\Delta_{g_p}e^{nH_p+\log f +n\log\frac{r_p^2}{r^2}-nu_p}(p)\\
		+&\mathcal{O}\left( s^3+\frac{1}{L^4}+\frac{1}{L^n}+\frac{s^2}{L^{n-2}} \right)\\
		=&C_0 e^{nH_p+n\mathcal{E}+\log f }(p)\left( 1+ \frac{C_2s^2}{C_0}\left( \Delta_{g_p}\left( H_p+\frac{1}{n}\log f \right)+n\left| \nabla_{g_p}\left( H_p+\frac{1}{n}\log f \right)  \right|^2 -\frac{R_g}{6(n-1)}  \right) (p) \right)\\
		+&\mathcal{O}\left( s^3+\frac{1}{L^4}+\frac{1}{L^n}+\frac{s^2}{L^{n-2}} +\frac{s^2}{L^2}\right).
	\end{split}
\]

For $n=2$, we have
\[
	\begin{split}
	&\lambda^{-2n}s^n\int_{B_{Ls}(p)} fe^{n\varphi_{p,s}}d\mu\\
		=&C_0 e^{nH_p+n\mathcal{E}+\log f }(p)\left( 1+ \frac{2s^2\log L}{\lambda^2}\left( \Delta_{g_p}\left( H_p+\frac{1}{n}\log f \right)+n\left| \nabla_{g_p}\left( H_p+\frac{1}{n}\log f \right)  \right|^2 -\frac{R_g}{6(n-1)}  \right)(p)  \right)\\
		+&\mathcal{O}\left( s^2+\frac{1}{L^2} + \frac{s^2\log L}{L^2}\right).
	\end{split}
\]
We set $L=s^{-\frac{1}{2}}\left( \log\frac{1}{s} \right)^{\frac{1}{4}}$ when $n>2$, so
\[
.\]
\[
\begin{split}
	E[\varphi_{p,s}]=&
	-(n-1)!\gamma_n\log\frac{\gamma_n f(p)}{2^n}-\frac{n!\gamma_n}{2}\left( H_p(p)+\frac{C_{\log}}{C_0} \right) \\
	-&(n-1)!\gamma_n\frac{C_2s^2}{C_0}\left( \Delta\left( H_p+\frac{1}{n}\log f \right)+n\left| \nabla\left( H_p+\frac{1}{n}\log f \right)  \right|^2 -\frac{R_g}{6(n-1)}  \right)(p) +\mathcal{O}\left( s^2\left( \log\frac{1}{s} \right) ^{-1} \right) .
\end{split}
\]

On the other hand,  we set $L=s^{-1}\left( \log\frac{1}{s} \right)^{-\frac{1}{4}}$ in the $n=2$ case ,  then we have
\[
\begin{split}
	E[\varphi_{p,s}]=&
	-4\pi\log f(p)-4\pi\left( H_p(p)+1+\log\pi \right) \\
	-&\frac{32\pi}{lf(p)}s^2\log\frac{1}{s}\left( \Delta\left( H_p+\frac{1}{2}\log f \right)+2\left| \nabla\left( H_p+\frac{1}{2}\log f \right)  \right|^2 -\frac{R_g}{6}  \right)(p) +\mathcal{O}\left( s^2\left( \log\frac{1}{s} \right) ^{\frac{1}{2}} \right)\\
=&-4\pi\log f(p)-4\pi\left( H_p(p)+1+\log\pi \right) \\
-&\frac{16\pi}{lf(p)}s^2\log\frac{1}{s}\left( \Delta\log f +4\left| \nabla\left( H_p+\frac{1}{2}\log f \right)  \right|^2 +2Q-R_g  \right)(p) +\mathcal{O}\left( s^2\left( \log\frac{1}{s} \right) ^{\frac{1}{2}} \right).
\end{split}
\]
In the last equality we have used the following identity for the two dimensional case (see \cite{DJLW1997}).
\begin{equation}\label{2dim}
	\Delta H_p(p) + \frac{R_g(p)}{3}=Q(p).
\end{equation}
This can be deduced directly by using
\[
	\int_{\partial B_r(p)}\frac{\partial G_p}{\partial n}d\sigma = -4\pi+\int_{B_r(p)}Qd\mu
,\]
and
\[
	\int_{\partial B_r(p)}\frac{\partial G_p}{\partial n}d\sigma = \int_{\partial B_r(p)}\left( -\frac{2}{r} +\frac{\partial H_p}{\partial n}  \right) d\sigma=-\frac{2\Vol (\partial B_r(p))}{r}+\int_{B_r(p)}\Delta H_p d\mu
,\]
and
\[
	\Vol (\partial B_r(p))=2\pi r-\frac{R_g(p)\pi}{6}r^3+\mathcal{O}(r^4)
.\]
From now on, we fix $L=s^{-\frac{1}{2}}\left( \log\frac{1}{s} \right)^{\frac{1}{4}}$ for $n>2$ and $L=s^{-1}\left( \log\frac{1}{s} \right)^{-\frac{1}{4}}$ for $n=2$.

With the help of $\varphi_{p,s}$, we can now describe the global behavior of $u$ as $\inf_{x\in M} s(x)$ tend to zero.
\begin{theorem}\label{global app}
	Assume $u$ is a solution of the equation (\ref{sequence equation}) and $P$ has trivial kernel.
    Then for $E_\alpha$ and $s_{min}=\inf_{x\in M}s(x)$ are both small enough (that is, smaller than a given small positive number), there exist $s>0$ and  $p\in M$ such that,
    $$ \|u-\bar{u}-\varphi_{p,s} \|_{ W^{ \frac{n}{2},2 } (M) } + \|u-\bar{u}-\varphi_{p,s} \|_{ L^{ \infty } (M) }  \leq J(E_\alpha, s_{min})  ,$$
$J$ is positive number satisfying
\[
	\lim_{ E_\alpha, s_{min}\to 0 } J(E_\alpha, s_{min})=0
.\]
\end{theorem}
{\bf Proof. }
Combining the equation (\ref{sequence equation}) and (\ref{app phi}), we have that
\[
\begin{split}
	P\left( u-\varphi_{p,s} \right)=&\left( lf+g \right) e^{nu}- \frac{lf(p)s^n\chi_{B_{Ls}^{g_p}}e^{nu_p}}{(s^2+\lambda^2r_p^2)^n}+\mathcal{O}\left( \frac{\chi_{B_{Ls}^{g_p}}}{L^{n+2}s^n} + \frac{s^2\chi_{B_{Ls}^{g_p}}}{r_p^{n-1}}\right)  \\
\end{split}
.\]
Then by using Corollary \ref{stdbubble convergence}, we have
\begin{equation}\label{L1 conv}
	\int_M \left| P\left( u-\varphi_{p,s} \right)  \right|d\mu\leq J_1(E_\alpha, s_{min}),
\end{equation}
where $J_1$ is a positive number such that
\[
	\lim_{ E_\alpha, s_{min}\to 0 } J_1(E_\alpha, s_{min})=0
.\]

Since $P$ has trivial kernel, invoking the green representation formula, we have the following estimate for $x\in B_\delta(p)$ ($\delta$ is chosen as in Proposition \ref{decay on neck}).
\[
\begin{split}
	&\left( u-\varphi_{p,s} \right)(x)-\overline{\left( u-\varphi_{p,s} \right)}\\
	= &\int_MG_p(x,y)P\left( u-\varphi_{p,s} \right)(y)d\mu_y\\
	=&\int_{B_\delta(p)}2\log\frac{1}{|x-y|}P\left( u-\varphi_{p,s} \right)(y)d\mu_y + \mathcal{O}(J_1(E_\alpha, s_{min}))\\
	=&\int_{B_\delta(p)}2\log\frac{1}{|x-y|}\left( \left( lf+g \right) e^{nu}- \frac{lf(p)s^n\chi_{B_{2Ls}}}{(s^2+\lambda^2r^2)^n}(y)\right)d\mu_y  + \mathcal{O}(J_1(E_\alpha, s_{min}))\\
	+&\mathcal{O}\left( \frac{1}{L^2}+Ls^3+L^2s^2 \right)\\
	=&\int_{B_{\delta}(p)}2\log\frac{s}{|x-y|}P\left( u-\varphi_{p,s} \right)(y)d\mu_y+2\log\frac{1}{s}\int_{B_{\delta}(p)}P\left( u-\varphi_{p,s} \right)(y)d\mu_y+\mathcal{O}(J_1(E_\alpha, s_{min})).\\
\end{split}
\]

To compute the first term in the right hand of above formula, we use the estimate of $\hat{u}_{p,s}$. However, the $\hat{u}_{p,s}$ is define by using $exp_p^{g_p}$ instead of $exp_p$  here.
\[
\begin{split}
	&\left| \int_{B_{\delta}(p)}2\log\frac{s}{|x-y|}P\left( u-\varphi_{p,s} \right)(y)d\mu_y \right| \\
	\lesssim& \int_{B_{\frac{\delta}{s}}(0)}\left| \log |x-y|\left( lfe^{n\hat{u}_{p,s}(y)}-\frac{lf(p)\chi_{B_{L}(0)}(y)}{(1+\lambda^2 |y|^2)^n} \right)  \right|dy \\
	+& \int_{B_{\frac{\delta}{s}}(0)}\left| \log |x-y| \right|\left| \hat{h}_{p,s}(y)e^{n\hat{u}_{p,s}(y)} \right| dy + \mathcal{O}\left( \frac{1}{L^2}\log\frac{1}{s}+Ls^3\log\frac{1}{s} \right) \\
	\lesssim& \left\Vert lfe^{n\hat{u}_{p,s}(y)}-\frac{lf(p)\chi_{B_{L}(0)}(y)}{(1+\lambda^2 |y|^2)^n}\right\Vert_{L^1\left( B_{\frac{\delta}{s}}(0) \right)  }+\left\Vert \hat{h}_{p,s}(y)e^{n\hat{u}_{p,s}(y)}\right\Vert_{L^1\left( B_{\frac{\delta}{s}}(0) \right)  }\\
	+& \left\Vert lfe^{n\hat{u}_{p,s}(y)}-\frac{lf(p)\chi_{B_{L}(0)}(y)}{(1+\lambda^2 |y|^2)^n}\right\Vert_{L^\alpha\left( B_{\frac{\delta}{s}}(0) \right)  }+\left\Vert \hat{h}_{p,s}(y)e^{n\hat{u}_{p,s}(y)}\right\Vert_{L^\alpha\left( B_{\frac{\delta}{s}}(0) \right)  }\\
	+&\mathcal{O}\left( \frac{1}{L^2}\log\frac{1}{s}+Ls^3\log\frac{1}{s}\right).\\
\end{split}
\]
One can see the  $L_1$ terms tend to $0$ as $E_\alpha$ and $s_{min}$ tend to $0$ by using Corollary \ref{stdbubble convergence}. The $L^\alpha$ terms also tend to $0$ once we know the $L^1$ terms tend to $0$, because we have the global $L^\infty$ control by Corollary \ref{global upp}.

It is remained to estimate the second term. By Corollary \ref{out vol},
\[
\begin{split}
	&2\log\frac{1}{s}\left| \int_{B_{\delta}(p)}P\left( u-\varphi_{p,s} \right)(y)d\mu_y \right| \\
	=&2\log\frac{1}{s}\left| \int_{M\setminus B_{\delta}(p)}P\left( u-\varphi_{p,s} \right)(y)d\mu_y \right| \\
	=&2\log\frac{1}{s}\left| \int_{M\setminus B_{\delta}(p)}\left( lf+h \right)e^{nu}(y)d\mu_y \right| \\
	\lesssim&s^{n(\nu-1)}\log\frac{1}{s}+\left( \int_{M\setminus B_\delta(p)} h^2e^{nu}d\mu\right)^{\frac{1}{2}} \left( \int_{M\setminus B_\delta(p)} e^{nu}d\mu\right)^{\frac{1}{2}}\\
	\lesssim &s^{\frac{n(\nu-1)}{2}}\log\frac{1}{s}.
\end{split}
\]
Combing above computations together, we have
\[
	\lim_{E_\alpha, s_{min}\to 0}\left\Vert\left( u-\varphi_{p,s} \right)-\overline{\left( u-\varphi_{p,s} \right)}\right\Vert_{L^\infty\left( B_\delta(p) \right)}=0
.\]
On the other hand,  we know from Proposition \ref{outside convergence}
\[
	\lim_{E_\alpha, s_{min}\to 0}\left\Vert\left( u-\varphi_{p,s} \right)-\overline{\left( u-\varphi_{p,s} \right)}\right\Vert_{L^\infty\left(M\setminus B_{\frac{\delta(p)}{2}} \right)}=0
.\]
Then we have
\[
	\lim_{s\to 0}\bar{\varphi}_{p,s}=0
.\]
Thus
\[
	\lim_{E_\alpha, s_{min}\to 0}\left\Vert u-\bar{u}-\varphi_{p,s} \right\Vert_{L^\infty\left( M \right)}=0
.\]
From the property of the GJMS operator, we have
\[
	\left\Vert u-\bar{u}-\varphi_{p,s} \right\Vert_{W^{\frac{n}{2},2}\left( M \right)}\leq C\left( \int_M \left( u-\bar{u}-\varphi_{p,s} \right)P\left( u-\varphi_{p,s} \right)d\mu +  \left\Vert u-\bar{u}-\varphi_{p,s} \right\Vert_{L^\infty\left( M \right)}  \right)
.\]
Then by the inequality (\ref{L1 conv}),
\[
	\lim_{E_\alpha, s_{min}\to 0}\left\Vert u-\bar{u}-\varphi_{p,s} \right\Vert_{W^{\frac{n}{2},2}\left( M \right)}\leq\lim_{E_\alpha, s_{min}\to 0}\left\Vert u-\bar{u}-\varphi_{p,s} \right\Vert_{L^\infty\left( M \right)}=0
.\]
\hfill $\Box$ \\
With the help of Theorem \ref{global app}, we can compute the value of $E[u]$.
\begin{corollary}\label{energy lim}
	Assume $u$ is a solution of the equation (\ref{sequence equation}) and $P$ has trivial kernel.
    Then for $E_\alpha$ and $s_{min}=\inf_{x\in M}s(x)$ are both small enough (that is, smaller than a given small positive number), there exist $s>0$ and  $p\in M$ such that,
    \[
	    E[u]=E[\varphi_{p,s}]+J(E_\alpha, s_{min})=\Lambda \left( g, f, p \right)+ J_1(E_\alpha, s_{min})
    .\]
Here $J$ and $J_1$ are positive numbers satisfying
\[
	\lim_{ E_\alpha, s_{min}\to 0 } J(E_\alpha, s_{min})=\lim_{ E_\alpha, s_{min}\to 0 } J_1(E_\alpha, s_{min})=0
.\]
And
\[
	\Lambda\left( g, f, p \right) = -(n-1)!\gamma_n\log \frac{\gamma_nf(p)}{2^n}-\frac{n!\gamma_n}{2}\left( H_p(p)+\frac{C_{\log}}{C_0} \right)
.\]
\end{corollary}
{\bf Proof. }
A direct calculation shows
\[
\begin{split}
	&\left| E[u]-E[\varphi_{p,s}]\right|\\
	=&\left| E[u-\bar{u}]-E[\varphi_{p,s}] \right|\\
	\leq& C\left(  \|u-\bar{u}-\varphi_{p,s} \|_{ W^{ \frac{n}{2},2 } (M) } + \|u-\bar{u}-\varphi_{p,s} \|_{ L^{ \infty } (M) }\right).
\end{split}
\]
This corollary follows from Theorem \ref{global app} and the explicit form of $E[\varphi_{p,s}]$.
\hfill $\Box$ \\

\section{Proof of Theorem \ref{two alternatives} and \ref{convergence theorem}}
\setcounter{equation}{0}

\medskip

In this section, we consider the convergence of the flow equation (\ref{parabolic equation}).

We begin with the Adams-Fontana inequality (see \cite{BCY1992}, \cite{F1993})
\begin{proposition}\label{AF Proposition}
	Let $(M, g)$ be a compact $n$ dimensional manifold whose GJMS operator $P$ is weakly positive with trivial kernel. Then for any $u\in W^{\frac{n}{2},2} (M)$, we have
  \begin{equation}\label{Adams-Fontana inequality}
	  \log\int_M e^{nu}d\mu\leq \frac{n}{2(n-1)!\gamma_n}\int_MuPud\mu+n\int_Mud\mu+C
  \end{equation}
  for some constant $C=C(M, g, n)$.
\end{proposition}

From the Adams-Fontana inequality, we know  that the functional $E$ has  uniform lower bound. More precisely, let $h$ be the solution of $$\begin{cases}Ph+Q=(n-1)!\gamma_n\\ \int_Mhd\mu=0,\end{cases}$$ and denote $u=h+v$. Then we have:
\begin{equation}
  \begin{split}
    E[u]&=\frac{n}{2}\int_MuPud\mu+n\int_M Q ud\mu-\bar{Q}\log\left(\int_Mfe^{nu}d\mu\right)\\
	&=\frac{n}{2}\int_MvPvd\mu+n\int_MuPhd\mu-\frac{n}{2}\int_MhPhd\mu\\
	&+n\int_M Q ud\mu-(n-1)!\gamma_n\log\left(\int_Mfe^{nh}e^{nv}d\mu\right)\\
	&\geq \frac{n}{2}\int_MvPvd\mu + n!\gamma_n\int_Mvd\mu- (n-1)!\gamma_n\log\left(\int_Mfe^{nv}d\mu\right)\\
	&-C\left(M,g,||f||_{L^{\infty}(M)},||h||_{W^{\frac{n}{2},2}{(M)}}\right)\\
	&\geq -L\left(M,g,||f||_{L^{\infty}(M)},||h||_{W^{\frac{n}{2},2} {(M)}}\right).
  \end{split}
\end{equation}
Now we recall the following result of  A. Fardoun and R. Regbaoui (Proposition 5.2 in \cite{FR2018}):
\begin{proposition}\label{FR local bound}
	Let $u\in C^{\infty}(M\times [0,T))$ be the solution of the equation (\ref{parabolic equation}) defined on a maximal interval $[0,T)$, and assume that the GJMS operator $P$ has trivial kernel. Then for any $L>0$ and any $T_0\in[0,T)$, there exists a positive constant $C$ depending on $L$, $g$, $f$, $T_0$ and $||u_0||_{W^{\frac{n}{2},2} (M)}$ such that,  if $$\inf_{t\in[0,T_0]}E[u(t)]\geq -L,$$
  then we have
  \begin{equation}\label{local time bound}
	  \sup_{t\in[0,T_0]}||u(t)||_{W^{\frac{n}{2},2} (M)}\leq C.
  \end{equation}
\end{proposition}
\begin{remark}
  A. Fardoun and R. Regbaoui only handled the case $f=constant$ and $n=4$. However, when $f$ is not a constant and $n$ is arbitrary, the proof is almost the same. For readers' convenience, we shall include the proof in Section $6$.
\end{remark}

Combining Proposition \ref{FR local bound} and part $(a)$ of Theorem \ref{brendle main}, we obtain the long time existence part in Theorem \ref{two alternatives}.

By a direct calculation, we know that
$$\frac{d}{dt}E[u(t)]=n\int_M \left(\frac{\partial u}{\partial t}\right)^2e^{nu(t)}d\mu.$$
Integrating this formula, we get
$$n\int^T_0\int_M\left(\frac{\partial u}{\partial t}\right)^2e^{nu(t)}d\mu dt=E(u_0)-E(u(T))\leq E[u_0]+ L<+\infty.$$
Since the right hand of the above formula is independent of $T$, we can conclude that
$$\int^{\infty}_0\int_M\left(\frac{\partial u}{\partial t}\right)^2e^{nu(t)}d\mu dt< \infty.$$
As a consequence, we have
\begin{equation}\label{remiander term small}
  \liminf_{t\to+\infty}\int_M\left(\frac{\partial u}{\partial t}\right)^2(t)e^{nu(t)}d\mu = 0.
\end{equation}
We want to show
\begin{proposition}\label{lim noinf}
	Assume $u(t)$ is a solution of the flow equation (\ref{parabolic equation}), and the GJMS operator $P$ is weakly positive with trivial kernel, then
\begin{equation}\label{term small}
  \lim_{t\to+\infty}\int_M\left(\frac{\partial u}{\partial t}\right)^2(t)e^{nu(t)}d\mu = 0.
\end{equation}
\end{proposition}
We need the following lemma.
\begin{lemma}
	Assume $u(t)$ is a solution of the flow equation (\ref{parabolic equation}), then there exits $C=C(u(0), M, g)$ such that
	\[
		\left( \int_M v^6d\mu_{g_t} \right) ^{\frac{1}{3}}\leq C\left( \int_M vP_{g_t}v d\mu_{g_t}+ \left( \int_M Q_{g_t}vd\mu_{g_t} \right)^2 \right)
	,\]
for any $v\in W^{\frac{n}{2},2}(M)$, where $g_t=e^{2u(t)}g$, $P_{g_t}$ and $Q_{g_t}$ is the GJMS operator and $Q$ curvature associated to $g_t$.
\end{lemma}
{\bf Proof. }
By a direct calculation, we have
\[
	\begin{split}
		E_{g_t}[v]&=\frac{n}{2}\int_MvP_{g_t}vd\mu_{g_t}+n\int_M Q_{g_t} vd\mu_{g_t}-(n-1)!\gamma_n\log\left(\int_Mfe^{nv}d\mu_{g_t}\right)\\
			  &=\frac{n}{2}\int_MvPvd\mu + n\int_M\left( Pu(t)+Q \right)vd\mu-(n-1)!\gamma_n\log\left( \int_M fe^{n(u(t)+v)}d\mu \right)\\
			  &=\frac{n}{2}\int_M(u(t)+v)P(u(t)+v)d\mu+n\int_MQ(u(t)+v)d\mu-(n-1)!\gamma_n\log\left( \int_Mfe^{n(u(t)+v)}d\mu \right)\\
			  &-\frac{n}{2}\int_Mu(t)Pu(t)d\mu-n\int Qu(t)d\mu\\
			  &=E[u(t)+v]-E[u(t)]-(n-1)!\gamma_n\log\left( \int_Mfe^{nu(t)}d\mu \right).
	\end{split}
\]
Since $u(t)$ is a solution of the flow equation, we have a uniform upper bound of \[
	E[u(t)]+(n-1)!\gamma_n\log\left(\int_M fe^{nu(t)} d\mu \right)
,\]
independent of $t$. So there exists a uniform positive constant $C$, such that
\[
	(n-1)!\gamma_n\log\left( \int_MCfe^{nv}d\mu_{g_t} \right)\leq \frac{n}{2}\int_MvP_{g_t}vd\mu_{g_t}+n\int_M Q_{g_t}vd\mu_{g_t}
.\]
Then if we assume $\int_M wP_{g_t}w d\mu_{g_t}=1$ and $\int_M Q_{g_t}wd\mu_{g_t}=0$,
we have \[
	\int_M e^{nw}d\mu_{g_t}+\int_M e^{-nw}d\mu_{g_t}\leq C(n, f, u(0),M)
.\]
From element calculus, we know
\[
	|w|^6\leq C_1 + C_2(e^{nw}+e^{-nw})
,\]
for some $C_1$, $C_2$ depending only on $n$.
Thus we can conclude that
\begin{equation}\label{normlized estimate}
	\int_M|w|^6d\mu_{g_t}\leq C(n, f, u(0), M).
\end{equation}
Then for general $v\in W^{\frac{n}{2},2}(M)$, if
\[
	\int_M vP_{g_t}v d\mu_{g_t}=0
,\]
from the trivial kernel assumption, we know $v$ is a constant function. So we have
\[
	\left(  \int_M v^6d\mu_{g_t}\right) ^{\frac{1}{3}} = Vol(M, g)^{\frac{1}{3}}v^2= \frac{Vol(M, g)^{\frac{1}{3}}}{\left( (n-1)!\gamma_n \right) ^2}\left( \int_M Qvd\mu_{g_t} \right) ^2
.\]

If
\[
	\int_M vP_{g_t}v d\mu_{g_t}\neq 0
,\]
we let
\[
	w=\left( \int_M vP_{g_t}v d\mu_{g_t} \right) ^{-\frac{1}{2}}\left( v-\frac{\int_M Q_{g_t}vd\mu_{g_t}}{\int_M Q_{g_t}d\mu_{g_t}} \right)
.\]
Apply (\ref{normlized estimate}) to this $w$, we have
\[
	\left( \int_M \left( v-\frac{\int_M Q_{g_t}vd\mu_{g_t}}{\int_M Q_{g_t}d\mu_{g_t}} \right)^6d\mu_{g_t}  \right)^{\frac{1}{3}}\leq C\left( \int_M vP_{g_t}v d\mu_{g_t} \right)^2
.\]
This lemma then follows from the triangle inequality.

\hfill $\Box$ \\

{\bf Proof of Proposition \ref{lim noinf}:}
A direct calculation shows
\[
\begin{split}
	\frac{\partial}{\partial t}\left( Q_{g_t}-l_tf \right) =&-P_{g_t}(Q_{g_t}-l_tf)+nQ_{g_t}(Q_{g_t}-l_tf)+nl_tf\int_M\frac{f}{\bar{f}}(Q_{g_t}-l_tf)d\mu_{g_t}
\end{split}
.\]
where
\[
	l_t=\frac{\int_M Q_gd\mu_{g_t}}{\int_M fd\mu_{g_t}}=\frac{(n-1)!\gamma_n}{\int_Mfd\mu_{g_t}}
.\]
So
\begin{equation}
\begin{split}
	\frac{\partial}{\partial t}\int_M (Q_{g_t}-l_tf)^2 d\mu_{g_t}=&-2\int_M(Q_{g_t}-l_tf)P_{g_t}(Q_{g_t}-l_tf)d\mu_{g_t}+n\int_M(Q_{g_t}-l_tf)^3d\mu_{g_t}\\
	&+2n\int_M l_tf(Q_{g_t}-l_tf)^2d\mu_{g_t}+2n!\gamma_n\left( \int_M \frac{f}{\bar{f}} (Q_{g_t}-l_tf)d\mu_{g_t} \right)^2\\
	\leq& C\left( \int_M (Q_{g_t}-l_tf)^2d\mu_{g_t}+ \left( \int_M(Q_{g_t}-l_tf )^2d\mu_{g_t}\right) ^3 \right) .
\end{split}
\end{equation}
Here we have used the following simple inequality:
\[
	\int_Mw^3d\mu\leq\left( \int_Mw^6 d\mu \right) ^{\frac{1}{4}}\left( \int_Mw^2 d\mu\right) ^{\frac{3}{4}}\leq \varepsilon\int_Mw^6d\mu+C(\varepsilon)\left( \int_Mw^2d\mu \right) ^3
.\]
Denote $y(t)=\int_M (Q_{g_t}-l_tf)^2d\mu_{g_t}=\int_M\left( \frac{\partial u}{\partial t} \right)^2(t)d\mu_{g_t} $, then we have
\[
	\frac{dy}{dt}\leq Cy(1+y^2)\leq Cy(1+y)^2
.\]
Integrating it, we obtain
\[
	C\int_{t_1}^{t_2}y(t)dt\geq \frac{1}{1+y(t_1)}-\frac{1}{1+y(t_2)}
.\]
Thus
\[
	\liminf_{t_2\to\infty}\frac{1}{1+y(t_2)}\geq \sup_{t\geq 0}\left( \frac{1}{1+y(t_1)}+ C\int_{t_1}^\infty y(t)dt \right)\geq \limsup_{t\to\infty} \left( \frac{1}{1+y(t_1)}+ C\int_{t_1}^\infty y(t)dt \right)=1
.\]
So
\[
	\lim_{t\to\infty}y(t)=\limsup_{t\to\infty}y(t)=0
.\]
\hfill $\Box$ \\

We denote the volume concentration radius function of $u(t)$ by $s_t$ in order to state our main result of this section.
\begin{theorem}\label{aaa}
	Assume $u(t)$ is a solution of the flow equation (\ref{parabolic equation}), if 	
	\[
	\limsup_{t\to\infty}\inf_{x\in M} s_t(x)>0
	,\]
	then there exists a smooth function $u_{\infty}$ satisfying the equation (\ref{elliptic equation}) such that
	\begin{equation}\label{strong ss}
		\lim_{t\to\infty}\|u(t)-u_\infty\|_{W^{n,2}(M)}=0
	\end{equation}
\end{theorem}
{\bf Proof. }
Note that
	\[
	\limsup_{t\to\infty}\inf_{x\in M} s_t(x)>0
	,\]
	the argument below Definition \ref{radius} allows us to find a sequence $\{t_k\}_{k\in \mathbb{N}}$, such that
	\[
		\|u(t_k)\|_{W^{n,2}(M)}\leq C
	,\]
	for some constant $C$ independent of $k$. Then there exists a function $u_\infty\in W^{n,2}(M)$, such that after passing to a subsequence, which we still denote by $\{t_k\}_{k\in \mathbb{N}}$:
\[
	\lim_{k\to\infty}\left(  \|u(t_k)-u_\infty\|_{L^\infty(M)}+\|u(t_k)-u_\infty\|_{W^{\frac{n}{2},2}(M)}\right) =0
.\]
Then we  argue as L. Sun et al.  \cite{SZ2021} by using the \L ojiasiewicz-Simon inequality.
Actually, the general result proved in \cite{FM2019} says:

There are positive constants $\beta>0$ and $\theta\in[1/2,1)$ such that
\[
	\left| E[u(t)]-E[u_\infty] \right|^\theta\leq \left(  \int_M \left( \frac{\partial u}{\partial t} \right)^2(t)e^{2nu(t)}d\mu
\right)^{\frac{1}{2}}
.\]
for $\|u(t)-u_\infty\|_{W^{\frac{n}{2},2}(M)}\leq \beta$.

On the other hand, using Trudinger's inequality, we have
\[
	\int_M e^{nqu(t)}d\mu\leq C(n, q,\beta)
.\]
for all $q>1$. Thus $u_t$ has a uniform lower bound as long as $\|u(t)-u_\infty\|_{W^{\frac{n}{2},2}(M)}\leq \beta$. Combing this fact with Proposition \ref{lim noinf}, we have
\[
	\|u(t)\|_{L^\infty(M)}\leq\|u(t)\|_{W^{n,2}(M)}\leq C(n,q,\beta)
,\]
for some $C(n,q,\beta)$ independent of $t$.
So we have
\[
	\left| E[u(t)]-E[u_\infty] \right|^\theta\leq \left(  \int_M \left( \frac{\partial u}{\partial t} \right)^2(t)e^{nu(t)}d\mu
\right)^{\frac{1}{2}}
.\]
for $\|u(t)-u_\infty\|_{W^{\frac{n}{2},2}(M)}\leq \beta$.

A direct calculation shows
\[
\begin{split}
	-\frac{d}{dt}(E[u(t)]-E[u_\infty])^{1-\theta}
	=&-(1-\theta)(E[u(t)]-E[u_\infty])^{-\theta}\frac{d}{dt}E[u(t)]\\
	=&(1-\theta)n(E[u(t)]-E[u_\infty])^{-\theta}\int_M \left( \frac{\partial u}{\partial t} \right)^2(t)e^{nu(t)}d\mu \\
	\geq &C\left( \int_M \left( \frac{\partial u}{\partial t} \right)^2(t)e^{nu(t)}d\mu  \right)^{\frac{1}{2}} \geq \left( \int_M \left( \frac{\partial u}{\partial t} \right)^2(t)d\mu  \right)^{\frac{1}{2}}.
\end{split}
\]
Thus
\[
	\frac{d}{dt}\|u(t)-u_\infty\|_{L^2(M)}\leq C\left( \int_M \left( \frac{\partial u}{\partial t} \right)^2(t)d\mu  \right)^{\frac{1}{2}}
.\]
So we have
\[
\begin{split}
	&\|u(\tau)-u_\infty\|_{L^2(M)}\\
	\leq &\|u(t_k)-u_\infty\|_{L^2(M)}+\int_{t_k}^\tau\left( \int_M \left( \frac{\partial u}{\partial t} \right)^2(t)d\mu  \right)^{\frac{1}{2}}dt\\
	\leq &\|u(t_k)-u_\infty\|_{L^2(M)}+\left( E[u(t_k)]-E[u_\infty] \right)^{1-\theta},
\end{split}
\]
once $\|u(t)-u_\infty\|_{W^{\frac{n}{2},2}(M)}\leq \beta$ for all $t\in [t_k, \tau]$.

We define
\[
	\tau_k=\sup\{\tau:\|u(t)-u_\infty\|_{W^{\frac{n}{2},2}(M)}\leq \beta \text{ for all } t\in [t_k, \tau]\}
.\]
If $\tau_k<\infty$, we have
\[
	\|u(\tau_k)-u_\infty\|_{W^{\frac{n}{2},2}(M)}=\beta\qquad\lim_{k\to\infty} \|u(\tau_k)-u_\infty\|_{L^2(M)}=0
.\]
And this is contradict to $\|u(\tau_k)\|_{W^{n,2}(M)}\leq C(n,q,\beta)$.
So we have $s_k=\infty$ for $k$ large enough.

Thus there exist a number $\tau_0>0$ such that
\[
\|u(\tau)\|_{W^{n,2}(M)}\leq C(n,q,\beta)\qquad \lim_{\tau\to\infty} \|u(\tau)-u_\infty\|_{L^2(M)}=0
,\]
for all $\tau\in [\tau_0, \infty)$.
This implies
\[
	\lim_{t\to\infty}\left(  \|u(t)-u_\infty\|_{L^\infty(M)}+\|u(t)-u_\infty\|_{W^{\frac{n}{2},2}(M)}\right) =0
.\]
The formula (\ref{strong ss}) then follows from element elliptic estimate as in the proof of Proposition \ref{strong compactness}.
\hfill $\Box$ \\

Now we can give a proof of Theorem \ref{two alternatives}.

{\bf Proof of Theorem \ref{two alternatives}:}
If $\limsup_{t\to\infty}\inf_{x\in M}s_t(x)>0$, Theorem \ref{aaa} show the alternative (a) holds. If $\lim_{t\to\infty}\inf_{x\in M}s_t(x)=0$. Corollary \ref{energy lim} tells the alternative (b) holds.
\hfill $\Box$ \\

{\bf Proof of Theorem \ref{convergence theorem}:}
From the assumption of this theorem and the explicit form of $E[\varphi_{p,s}]$, we can find a positive number $s_0$ such that
\[
	E[\varphi_{p_0,s_0}]<\inf_{p\in M}\Lambda\left( g,f,p \right)
.\]
We now choose $\varphi_{p_0,s_0}$ as the initial data, the alternative (b) in Theorem \ref{two alternatives} cannot hold, and the flow converges.
\hfill $\Box$ \\
  \section{Proof of Theorem \ref{large convergence theorem}}

  \setcounter{equation}{0}

  \medskip
  In this section, we study  the flow equation (\ref{parabolic equation}) under the assumption
  \[
	  \lim_{t\to\infty}\inf_{x\in M}s_t(x)=0
  .\]
  By applying Theorem \ref{global app}, we know there exist $p_t$, $s_t$ such that
  \[
	  \lim_{t\to\infty}\left( \|u(t)-\bar{u}(t)-\varphi_{p_t,s_t} \|_{ W^{ \frac{n}{2},2 } (M) } + \|u(t)-\bar{u}(t)-\varphi_{p_t,s_t} \|_{ L^{ \infty } (M) } \right)=0
  .\]
  We also have $\lim_{t\to\infty}s_t=0$.

  We now choose $(p_t^*, s_t^*)\in M\times\mathbb{R}^+$ such that
  \begin{equation}\label{mini prob}
	  \int_M\left( u(t)-\varphi_{p_t^*,s_t^*} \right) P\left( u(t)-\varphi_{p_t^*,s_t^*} \right)d\mu\leq \inf_{(p_t^*, s_t^*)\in M\times\mathbb{R}^+}\int_M\left( u(t)-\varphi_{p,s}\right) P\left( u(t)-\varphi_{p,s} \right)d\mu.
  \end{equation}

  Thus we have
  \[
	  \lim_{t\to\infty}\int_M\left( \varphi_{p_t, s_t}-\varphi_{p_t^*,s_t^*} \right) P\left( \varphi_{p_t, s_t}-\varphi_{p_t^*,s_t^*} \right)d\mu =0
  .\]
Using this, we can check
\[
	\lim_{t\to\infty} \frac{d_g(p_t,p_t^*)}{s_t^*}=0\qquad\quad \lim_{t\to\infty}\frac{s_t}{s_t^*}=1
.\]
Then
  \[
	  \lim_{t\to\infty}\left( \|u(t)-\bar{u}(t)-\varphi_{p_t^*,s_t^*} \|_{ W^{ \frac{n}{2},2 } (M) } + \|u(t)-\bar{u}(t)-\varphi_{p_t^*,s_t^*} \|_{ L^{ \infty } (M) } \right)=0
  .\]

  The inequality also implies
  \[
	  \int_M\left( u(t)-\varphi_{p_t^*,s_t^*} \right)P \frac{\partial \varphi_{e_{p_t^*}(y),s_t^*}}{\partial y}\Bigg|_{y=0}d\mu=0
  ,\]
  and
  \[
	  \int_M\left( u(t)-\varphi_{p_t^*,s_t^*} \right)P \frac{\partial \varphi_{p_t^*,s_t^*+z}}{\partial z}\Bigg|_{s=0}d\mu=0
  .\]

  To proceed, we need the following Poincar\'e inequality by S. Brendle's blow up argument (see\cite{S2005}) .
  \begin{lemma}\label{poin}
Assume $u$ is a solution of the equation (\ref{sequence equation}).
    Then for $E_\alpha$ and $s_{min}=\inf_{x\in M}s(x)$ are both small enough (that is, smaller than a given positive real number), there exist a uniform constant $\beta>0$, such that
	\[
		\left( 1-\beta \right) \int_M   \omega P\omega\geq n!\gamma_n  \frac{\int_M \omega^2fe^{nu}d\mu}{\int_M fe^{nu}d\mu}
	,\]
for any $\omega\in W^{\frac{n}{2},2}\left( M \right)$ satisfying the following three constraints:
\[
	\int_M \omega f e^{nu}d\mu=0
.\]
    \[
		\int_M\omega P\frac{\partial \varphi_{e_p(y),s}}{\partial y}\Bigg|_{y=0}d\mu=0
.\]
And
\[
\int_M\omega P\frac{\partial \varphi_{p,s+z}}{\partial z}\Bigg|_{z=0}d\mu=0
.\]
\end{lemma}
{\bf Proof. }
Supposing this is not true. We can assume there exists a sequence of function $\{u_k\}_{k\in \mathbb{N}}$ satisfying:
\begin{equation}
  \begin{cases}
  Pu_k +Q = (l_kf_k+ h_k)e^{nu_k},
\\
\int_M l_kf_ke^{nu_k}d\mu = \int_M Q d\mu = (n-1)!\gamma_n
\\
\lim_{k\to \infty} E_{\alpha,k} =\lim_{k\to\infty} \int_M g_k^\alpha e^{nu_k}d\mu = 0,
  \end{cases}
\end{equation}
and a sequence of function $\{\omega_k\}_{k\in \mathbb{N}}$ satisfying the three constraints in the lemma such that:
\[
		\int_M   \omega_k P\omega_kd\mu\leq n!\gamma_n  \frac{\int_M \omega_k^2fe^{nu_k}d\mu}{\int_M fe^{nu_k}d\mu}=1
.\]
Assume $p_k$, $s_k$ is a pair of parameter of $u_k$ satisfying the conclusion of  Proposition \ref{stdbubble convergence}, and denote $\hat{u}_k=\widehat{\left( u_k \right) }_{p_k, s_k}$ and $\hat{\omega}_k=\widehat{\left( \omega_k \right) }_{p_k, s_k}$.
Then we have
\[
	\int_{B_{\frac{\delta}{s}}(0)}\hat{\omega}_k P_{g_{p,s}}\hat{\omega}_kd\mu_{p,k}\leq C\left( n, g, f \right)
,\]
and
\[
	\int_{B_{\frac{\delta}{s}}(0)}\hat{\omega}_k^2e^{n\hat{u}_k} d\mu_{p,k}\leq C\left( n, g, f \right)
.\]
Hence, if we take the weak limit as $k\to\infty$, we obtain a function $\hat{\omega}\in W^{\frac{n}{2},2}\left( \mathbb{R}^n \right)$ such that
\begin{equation}\label{inve poin}
	\int_{\mathbb{R}^n}\hat{\omega}\left( -\Delta_0 \right) ^{\frac{n}{2}}\hat{\omega}dx\leq n!\gamma_n\int_{\mathbb{R}^n}\frac{\hat{\omega}^2}{\left( 1+\lambda^2|x|^2 \right)^n}dx \Bigg/\int_{\mathbb{R}^n}\frac{1}{\left( 1+\lambda^2|x|^2 \right)^n}dx=1.
\end{equation}

Moreover, we have the following estimate from a direct calculation.
\[ \begin{split}
	&P\frac{\partial}{\partial y}\varphi_{e_p(y),s}(x)+\frac{2nlf(p)s^n\lambda^2e_p^{-1}(x)
 \cdot y}{(s^2+\lambda^2r_p^2(x))^{n+1}}\\
	=&\frac{\partial}{\partial y}\left( \left( e^{nu_{e_p(y)}}-1 \right)P_{g_{e_p(y)}}\varphi_{e_p(y), s}  \right) + \frac{\partial}{\partial y}\left( P_{g_{e_p(y)}}\left( \varphi_{e_p(y),s}+\log\left( \frac{s^2}{\lambda^2}+r_{e_p(y)}^2 \right) \right) \right) \\
+&\frac{\partial}{\partial y}\left(  \left( \left( -\Delta_{g_{e_p(y)}} \right)^{\frac{n}{2}}-P_{g_{e_p(y)}}  \right)\log\left( \frac{s^2}{\lambda^2}+r_{e_p(y)}^2 \right) \right)  \\
+&\frac{\partial}{\partial y}\frac{lf(p)s^n}{\left( s^2+\lambda^2r_{e_p(y)}^2 \right)^n}+\frac{2nlf(p)s^n\lambda^2e_p^{-1}(x) \cdot y }{(s^2+\lambda^2r_p^2(x))^{n+1}}+\mathcal{O}\left( \frac{1}{L^{2n+1}s^{n+1}}\right) \\
=& \mathcal{O}\left(  \frac{1}{r_p^{n-1}}+\frac{1}{L^{n+3}s^{n+1}}\right),
\end{split}\]
and
\[\begin{split}
	&P\frac{\partial}{\partial s}\varphi_{p,s}(x)+\frac{nlf(p)s^{n-1}\left( s^2 -\lambda^2r_p^2\right) }{(s^2+\lambda^2r_p^2(x))^{n+1}}\\
	=&\frac{\partial}{\partial s}\left( \left( e^{nu_{p}}-1 \right)P_{g_{p}}\varphi_{p, s}  \right) + \frac{\partial}{\partial s}\left( P_{g_{p}}\left( \varphi_{p,s}+\log\left( \frac{s^2}{\lambda^2}+r_{p}^2 \right) \right) \right) \\
+&\frac{\partial}{\partial s}\left(  \left( \left( -\Delta_{g_{p}} \right)^{\frac{n}{2}}-P_{g_p}  \right)\log\left( \frac{s^2}{\lambda^2}+r_{p}^2 \right) \right) +\mathcal{O}\left( \frac{1}{L^{2n+1}s^{n+1}} \right)  \\
=&\mathcal{O}\left(  \frac{1}{r_p^{n-1}}+\frac{1}{L^{n+2}s^{n+1}}\right).
\end{split}\]
Here we use $e_p(y)$ to denote $exp_p^{g_p}(y)$.

Thus, we have the following constraints hold:
\[
	\int_{\mathbb{R}^n}\frac{\hat{\omega}}{\left( 1+\lambda^2|x|^2 \right)^n}dx=0
,\]
and
\[
	\int_{\mathbb{R}^n}\frac{\left(  1-\lambda^2|x|^2\right) \hat{\omega}}{\left( 1+\lambda^2|x|^2 \right)^{n+1}}dx=0
,\]
and
\[
	\int_{\mathbb{R}^n}\frac{x\hat{\omega}}{\left( 1+\lambda^2|x|^2 \right)^{n+1}}dx=0
.\]
By the conformal covariance of the $GJMS$ operator, and the spectral property of the $GJMS$ operator on the round $S^n$, the equation (\ref{inve poin}) cannot hold. A contradiction.

\hfill $\Box$ \\
We denote that
\[
	\phi(t)=\varphi_{p_t^*,s_t^*} + B_t \qquad \quad \omega(t)=u(t)-\phi(t)
,\]
where
\[
	B_t=\frac{\int_M\left( u(t)-\varphi_{p_t^*,s_t^*}\right) fd\mu_{g_t} }{\int_Mfd\mu_{g_t}}
.\]
Thus
\[
	\int_M\omega(t)fd\mu_{g_t}=0
.\]
We also have
  \[
	  \int_M\omega(t)P \frac{\partial \varphi_{e_{p_t^*}(y),s_t^*}}{\partial y}\Bigg|_{y=0}d\mu=0
,\]
and
  \[
	  \int_M\omega(t)P \frac{\partial \varphi_{p_t^*,s_t^*+z}}{\partial z}\Bigg|_{z=0}d\mu=0
  .\]
Then we can calculate
\[
\begin{split}
	&E[u(t)]-E[\phi(t)]\\
	=&n\int_M \omega(t)\left( P u(t)+ Q\right)  d\mu - \frac{n}{2}\int_M \omega(t)P\omega(t)d\mu + (n-1)!\gamma_n\log \left( \frac{\int_Me^{-n\omega(t)}fd\mu_{g_t}}{\int_M fd\mu_{g_t}} \right) \\
	=&n\int_M \omega(t)\frac{\partial u}{\partial t}(t) d\mu_{g_t}-\frac{n}{2} \int_M   \omega(t)P\omega(t)+(n-1)\gamma_n\log \left( 1-\frac{\int_M \left( n\omega(t)+\frac{n^2}{2}\omega^2(t) \right) fd\mu_{g_t}}{\int_Mfd\mu_{g_t}} \right)\\
	+&(n-1)!\gamma_n\log \left( \frac{\int_M \left( e^{-n\omega(t)}-1+n\omega(t)-\frac{n^2}{2}\omega^2(t) \right)fd\mu_{g_t} }{\int_M fd\mu_{g_t}} \right) \\
	=&n\int_M \omega(t)\frac{\partial u}{\partial t}(t) d\mu_{g_t}-\frac{n}{2}\left( \int_M   \omega(t)P\omega(t)-n!\gamma_n  \frac{\int_M \omega^2(t)fd\mu_{g_t}}{\int_M fd\mu_{g_t}} \right)\\
	+&\mathcal{O}\left( \frac{\int_M \omega^3fd\mu_{g_t}}{\int_M fd\mu_{g_t}} \right)
	\leq  C\int_M \left( \frac{\partial u}{\partial t} \right) ^2(t)d\mu_{g_t}.
\end{split}
\]
Notice that in the last inequality, we have used Lemma \ref{poin}.

We now state the main result of this section.
\begin{theorem}\label{energy low conv}
	Assume $u(t)$ is a solution of the flow equation (\ref{parabolic equation}), and denote $E_\infty=\lim_{t\to\infty}E[u(t)]$, if 	
\[
\lim_{t\to\infty}\inf_{x\in M} s_t(x)=0
,\]
then there does not exist a positive number $T>0$, such that	\begin{equation}\label{ci}
		E[\phi(t)]\leq E_\infty
\end{equation}
for all $t\in [T,\infty)$.
\end{theorem}

{\bf Proof:}

Supposing this is not true. We can assume there exists a positive number $T>0$ such that (\ref{ci}) holds for all $t\in [T,\infty)$. Thus we have
\[
	E[u(t)]-E[u_\infty]\leq C\left(  \int_M \left( \frac{\partial u}{\partial t} \right) ^2(t)d\mu_{g_t}\right)
,\]
for all $t\in [T, \infty)$.
\[
	\int_{\tau}^\infty\int_M \left( \frac{\partial u}{\partial t}\right)^2(t)d\mu_{g_t}dt\leq -C\frac{d}{d\tau}\int_\tau^{\infty}\int_M \left( \frac{\partial u}{\partial t}\right)^2(t)d\mu_{g_t}dt
.\]
So we have
\[
	\int_{\tau}^\infty\int_M \left( \frac{\partial u}{\partial t}\right)^2(t)d\mu_{g_t}dt\leq A^2 e^{-2\beta\tau}
,\]
for some uniform constants $A>0$ and $\beta>0$.
Then
\[
	\int_\tau^\infty\left( \int_M \left( \frac{\partial u}{\partial t}\right)^2(t)d\mu_{g_t} \right)^{\frac{1}{2}}dt=\sum_{i=0}^\infty\int_{\tau+i}^{\tau+i+1}\left( \int_M \left( \frac{\partial u}{\partial t}\right)^2(t)d\mu_{g_t} \right)^{\frac{1}{2}}dt\leq\sum_{i=0}^\infty Ae^{-\beta(\tau+i)}= \frac{Ae^{-\beta\tau}}{1-e^{-\beta}}
.\]
Thus
\[
	\begin{split}
		\int_{B_r(p)}e^{nu(\tau_1)}d\mu =&\int_{B_r(p)}e^{nu(\tau)}d\mu+n\int_\tau^{\tau_1}\int_{B_r(p)}\frac{\partial u}{\partial t}e^{u(t)}dt\\
&\leq\int_{B_r(p)}e^{nu(\tau)}d\mu+\int_\tau^{\tau_1}\left( \int_M \left( \frac{\partial u}{\partial t}\right)^2(t)d\mu_{g_t} \right)^{\frac{1}{2}}dt\\
&\leq\int_{B_r(p)}e^{nu(\tau)}d\mu +\frac{Ae^{-\beta\tau}}{1-e^{-\beta}}.
	\end{split}
\]
And this is contradict to the assumption $\lim_{t\to\infty}\inf_{x\in M} s(x)=0$.

\hfill $\Box$ \\
Now we can give a proof of the Theorem \ref{large convergence theorem}.

{\bf Proof of Theorem \ref{large convergence theorem}}
By our assumption \[
\Lambda\left( g, f, p \right)\equiv \Lambda
,\]
for some constant $\Lambda$.
Then Corollary \ref{energy lim} tells us
\[
	E_\infty=\lim_{t\to\infty}E[u(t)]=\Lambda
.\]
By the explicit form of $E[\varphi_{p_t^*,s_t^*}]$ and our assumption we can show
\[
	E[\phi(t)]=E[\varphi_{p_t^*,s_t^*}]< \Lambda=E_\infty
,\]
for $t$ large enough. Thus we have
\[
	\limsup_{t\to\infty}\inf_{x\in M}s_t(x)>0
,\]
by using Theorem \ref{energy low conv}. The convergence then follows from Theorem \ref{aaa}.
\hfill $\Box$ \\
\begin{remark}
	When $(M, g)$ is the standard round sphere $(S^n, g_{round})$, and $f\equiv 1$,  we know the standard bubble function $\Phi_{p,s}$ satisfying
	\[
		E[\Phi_{p,s}]\equiv \Lambda
	,\]
	for some constant $\Lambda$. If $\lim_{t\to\infty}\inf_{x\in M}s_t(x)=0$, Corollary \ref{energy lim} tells us $E_\infty=\Lambda$. Then we can argue as before by using Theorem \ref{energy low conv} and Theorem \ref{aaa} to show convergence in this special case. In this way, we recover the results of S. Brendle (see \cite{B2006}) and P.T. Ho (see\cite{H2010}).

\end{remark}

  \section{Proof of Proposition \ref{FR local bound}}

  \setcounter{equation}{0}

  \medskip

  In this section we introduce the notation $|A|=\int_A d\mu$ and assume $\int_M e^{nu(t)}d\mu\equiv\int_Md\mu=1$ for simplicity.

  By a direct calculation, we have
  $$\frac{d}{dt}\int_M ue^{nu}d\mu=-n\int_MuPud\mu-n\int_MQd\mu+n\frac{\bar{Q}}{\bar{f}}\int_Mfue^{nu}d\mu.$$
  Denote  $U=ue^{nu}+\frac{1}{n}e^{-1}$, then $U\geq0$ from elementary calculus.
  By the assumption, we have
  $$\frac{n}{2}\int_M uP ud\mu + n\int_M Qud\mu-\left(\int_MQd\mu\right)\log\left(\int_Mfe^{nu}d\mu\right)\geq -L.$$
  Therefore,
  \begin{equation}
    \begin{split}
      \frac{d}{dt}\int_MUd\mu &\leq-\frac{n}{2}\int_MuPud\mu+n\frac{\bar{Q}}{\bar{f}}\int_M Ud\mu-\frac{\bar{Q}}{e\bar{f}}\int_Mfd\mu\\
			      &+\left(\int_MQd\mu\right)\log\left(\int_Mfe^{nu} d\mu\right)+L\\
			      &\leq A(M,g,f)\int_M Ud\mu+C(M,g,f ,L) .
    \end{split}
  \end{equation}
  Then we can deduce:
  $$\int_Mue^{nu}d\mu\leq C(M,g,T_0,f).$$

  Denote $\varphi(z)=z\log z$. Then $\varphi$ is convex on $(0,+\infty)$, and it satisfies that for each $\lambda>1$ and $z>0$
  $$z=\frac{\varphi(\lambda z)}{\varphi(\lambda)}-\frac{\varphi(z)}{\log\lambda},$$
  which implies, because $\varphi(z)\geq-e^{-1}$ for any $z>0$,
  $$z\leq\frac{\varphi(\lambda z)}{\varphi(\lambda)}+\frac{e^{-1}}{\log\lambda}.$$
  For $t\in[0, T_0]$, let $A_t\subset M$ defined by
  $$A_t={x\in \{M:u(x,t)\geq\alpha_0\}},$$
  where
  $$\alpha_0=\frac{1}{n}\log\left(\frac{1}{2}\right).$$
  We shall demonstrate $|A_t|\geq C>0$ for some $C$ independent of $t$. If $|A_t|\geq 1$, there is nothing to prove. So we assume $|A_t|<1$.
  From Jensen's inequality it follows
  $$\varphi\left(\frac{1}{|A_t|}\int_{A_t}e^{nu}d\mu\right)\leq\frac{1}{|A_t|}\int_{A_t}\varphi(e^{nu})d\mu\leq\frac{C_1}{|A_t|}.$$
  Choosing $\lambda=\frac{1}{|A_t|}$ and $z=\int_{A_t}e^{nu}d\mu$, we have
  \begin{equation}
    \begin{split}
      \int_{A_t}e^{nu}d\mu
&\leq\frac{|A_t|}{\log\frac{1}{|A_t|}}\varphi\left(\frac{1}{|A_t|}\int_{A_t}e^{nu}d\mu\right)+\frac{e^{-1}}{\log\frac{1}{|A_t|}}\\
&\leq(C_1+e^{-1})\frac{1}{\log\frac{1}{|A_t|}}\leq\frac{C_2}{\log\frac{1}{|A_t|}}.
    \end{split}
  \end{equation}
  On the other hand, we have
  $$1=\int_Me^{nu}d\mu=\int_{A_t}e^{nu}d\mu+\int_{M\setminus A_t}e^{nu}d\mu.$$
  Since $e^{nu}<e^{\alpha_0}=\frac{1}{2}$ on $M\setminus A_t$, we get
  $$\frac{1}{2}\leq\int_{A_t}e^{nu}d\mu,$$
  which implies $|A_t|\geq e^{-2C_2}=C_3$. So we have  lower bound of $|A_t|$ which is independent of $t$.

  By using the elementary inequality $z\leq e^z$, we have
  $$\int_{A_t} u(t)d\mu\leq\frac{1}{n}\int_{A_t} e^{nu}d\mu\leq\frac{1}{n}.$$
  Then by the definition of the set $A_t$, for any $t\in[0,T_0]$, we conclude
  \begin{equation}
    \begin{split}
      \left|\int_M u(t)d\mu\right|
&\leq \left|\int_{A_t}u(t)d\mu\right|+\left|\int_{M\setminus A_t}u(t)d\mu\right|\\
&\leq C_4+\left|\int_{M\setminus A_t}u(t)d\mu\right|.
    \end{split}
  \end{equation}
  Using the Cauchy-Schwarz inequality and the Young's inequality, one get
  $$\left|\int_{M\setminus A_t}u(t)d\mu\right|\leq|M\setminus A_t|^{\frac{1}{2}}||u||_{L^2(M)},$$
  and
  $$\left(\int_M u(t)d\mu\right)^2\leq(1+\varepsilon)|M\setminus A_t|\lVert u(t) \rVert ^2_{L^2(M)}+C_4(1+\varepsilon^{-1}).$$
  Now, by Poincar\'e inequality we have

  $$\lVert u(t)\rVert^2_{L^2(M)}\leq K\int_M u(t)P u(t)d\mu+\frac{1}{\lvert M\rvert }\left(\int_Mu(t)d\mu\right) ^2.$$
  Combining above two inequality, we can conclude
  $$\left(1-\frac{(1+\varepsilon)|M\setminus A_t|}{|M|}\right)\lVert u(t) \rVert^2_{L^2(M)}\leq K\int_M u(t)P u(t)d\mu+\frac{C_4}{|M|}(1+\varepsilon^{-1}).$$
  That is,
  $$(\lvert A_t\rvert -\varepsilon\lvert M\setminus A_t\rvert)\lVert u(t)\rVert^2_{L^2(M)}\leq K|M|\int_M u(t)P u(t)d\mu+C_4(1+\varepsilon^{-1}).$$
  Choosing $\varepsilon=\frac{C_3}{2|M|}$ and observing that $|M\setminus A_t|\leq|M|$, we obtain
  \begin{equation}\label{poncare without mean value term}
    \lVert u(t)\rVert^2_{L^2(M)}\leq C_5\left(\int_M u(t)P u(t)d\mu +1\right).
  \end{equation}
  Since the functional $E$ is decreasing along the flow, we have
  $$\frac{n}{2}\int_M u(t)P u(t)d\mu + n\int_M Qu(t)d\mu-\bar{Q}\log\left(\int_Mfe^{nu}d\mu\right)\leq E[u_0],$$
  hence
  \begin{equation}\label{inverse poincare from monocinitity}
    \int_M u(t)P u(t)d\mu\leq \frac{1}{2C_5}\lVert u(t)\rVert^2_{L^2(M)}+C_6.
  \end{equation}
  It is follows from (\ref{poncare without mean value term}) and (\ref{inverse poincare from monocinitity}) that
  \begin{equation}\label{bound in last section}
    \lVert u(t)\rVert_{L^2(M)}\leq C_7.
  \end{equation}
  Combining (\ref{inverse poincare from monocinitity}) and (\ref{bound in last section}), we get the  estimate (\ref{local time bound}):
  $$\lVert u(t)\rVert_{W^{\frac{n}{2},2}(M)}\leq C_8.$$

  \par

\end{document}